\def\version{}

\documentclass[reqno,twoside,11pt]{amsart}
\usepackage{cite}
\usepackage{amsmath}
\usepackage{amsfonts}
\usepackage{amssymb}
\usepackage{epsfig}
\usepackage{verbatim}
\usepackage{color}

\usepackage{mathpazo}

\DeclareFontFamily{OT1}{cmss}{} \DeclareFontShape{OT1}{cmss}{m}{n} {<5> <6> <7> <8> <9> <10> <11> <12> <13> <14.4> cmss10}{}
\DeclareMathAlphabet{\cmss}{OT1}{cmss}{m}{n}

\newtheoremstyle{thm}{1.5ex}{1.5ex}{\itshape\rmfamily}{} {\bfseries\rmfamily}{}{2ex}{}

\newtheoremstyle{def}{1.5ex}{1.5ex}{\rmfamily\sl}{} {\bfseries\rmfamily}{}{2ex}{}

\newtheoremstyle{rem}{1.3ex}{1.3ex}{\rmfamily}{} {\bfseries\rmfamily}{}{2ex}{}

\newtheoremstyle{ass}{1.5ex}{1.5ex}{\rmfamily\sl}{} {\bfseries\rmfamily}{}{2ex}{}

\newenvironment{proofsect}[1] {\vskip0.1cm\noindent{\rmfamily\itshape#1.}}{\qed\vspace{0.15cm}}

\theoremstyle{thm}
\newtheorem{theorem}{Theorem}[section]
\newtheorem{lemma}[theorem]{Lemma}
\newtheorem{proposition}[theorem]{Proposition}

\newtheorem*{Main Theorem}{Main Theorem.}
\newtheorem{corollary}[theorem]{Corollary}

\newtheorem{assumption}[theorem]{Assumption}

\theoremstyle{rem}
\newtheorem{remark}[theorem]{{Remark}}

\numberwithin{equation}{section}

\renewcommand{\section}{\secdef\sct\sect}
\newcommand{\sct}[2][default]{\refstepcounter{section}
\addcontentsline{toc}{section}
{{\tocsection {}{\thesection}{\!\!\!\!#1\dotfill}}{}}
\vspace{0.7cm}
\centerline{ 
\scshape\arabic{section}.\ #1} \nopagebreak \vspace{0.2cm}}
\newcommand{\sect}[1]{
\vspace{0.4cm} \centerline{\large\scshape\rmfamily #1}
\vspace{0.2cm}}

\renewcommand{\subsection}{\secdef\subsct\sbsect}
\newcommand{\subsct}[2][default]{\refstepcounter{subsection}
\addcontentsline{toc}{subsection}
{{\tocsection{\!\!}{\hspace{1.2em}\thesubsection}{\!\!\!\!#1\dotfill}}{}}
\nopagebreak\vspace{0.45\baselineskip} {\flushleft\bf
\thesection.\arabic{subsection}~\bf #1.~}
\\*[3mm]\noindent
\nopagebreak}
\newcommand{\sbsect}[1]{
\vspace{0.1cm}\noindent
\textbf{#1.~}\vspace{0.1cm}}

\renewcommand{\subsubsection}{%
\secdef \subsubsect\sbsbsect}
\newcommand{\subsubsect}[2][default]{%
\refstepcounter{subsubsection} 
\addcontentsline{toc}{subsubsection}{{\tocsection{\!\!}
{\hspace{3.05em}\thesubsubsection}{\!\!\!\!#1\dotfill}}{}}
\nopagebreak
\vspace{0.15\baselineskip} \nopagebreak {\flushleft\rmfamily
\itshape\arabic{section}.\arabic{subsection}.\arabic{subsubsection}
\ \rmfamily #1\/.}\ }
\newcommand{\sbsbsect}[1]{\vspace{0.1cm}\noindent
\rmfamily \itshape
\arabic{section}.\arabic{subsection}.\arabic{subsubsection} \
\sffamily #1\/.\ }

\renewcommand{\caption}[1]{%
\vglue0.5cm
\refstepcounter{figure}
\begin{center}
\begin{minipage}[c]{0.8\textwidth}\small {\sc Fig.~\thefigure\ }#1\end{minipage}
\end{center}
}

\newcommand{\diam}{\operatorname{diam}}

\newcommand{\textd}{\text{\rm d}\mkern0.5mu}

\newcommand{\texte}{\text{\rm  e}\mkern0.7mu}

\newcommand{\BB}{\mathcal B}

\newcommand{\DD}{\mathcal D}

\newcommand{\FF}{\mathcal F}

\newcommand{\HH}{\mathcal H}

\newcommand{\LL}{\mathcal L}

\newcommand{\NN}{\mathcal N}

\newcommand{\E}{\mathbb E}

\newcommand{\G}{\mathbb G}

\newcommand{\N}{\mathbb N}

\newcommand{\BbbP}{\mathbb P}
\newcommand{\Q}{\mathbb Q}
\newcommand{\R}{\mathbb R}

\newcommand{\Z}{\mathbb Z}

\newcommand{\twoeqref}[2]{(\ref{#1}--\ref{#2})}

\newcommand{\cc}{{\text{\rm c}}}

\def\myffrac#1#2 in #3{\raise 2.6pt\hbox{$#3 #1$}\mkern-1.5mu\raise 0.8pt\hbox{$#3/$}\mkern-1.1mu\lower 1.5pt\hbox{$#3 #2$}}

\newcommand{\wh}{\widehat}
\newcommand{\wt}{\widetilde}

\newcommand{\laweq}{\,\overset{\text{\rm law}}=\,}

\begin{document}

\title[1D dynamical conductance model \hfill \version\hfill]
{\large An invariance principle for one-dimensional random\\walks among dynamical random conductances}

\author[\hfill  \version \hfill M.~Biskup]
{Marek~Biskup}
\thanks{\hglue-4.5mm\fontsize{9.6}{9.6}\selectfont\copyright\,\textrm{2018}\ \ \textrm{M.~Biskup.
Reproduction, by any means, of the entire
article for non-commercial purposes is permitted without charge.\vspace{2mm}}}
\maketitle

\vspace{-5mm}
\centerline{\textit{
Department of Mathematics, UCLA, Los Angeles, California, USA}}
\centerline{\textit{
Center for Theoretical Study, Charles University, Prague, Czech Republic}}


\vskip0.5cm
\begin{quote}
\footnotesize \textbf{Abstract:}
We study variable-speed random walks on~$\Z$ driven by a family of nearest-neighbor time-dependent random conductances $\{a_t(x,x+1)\colon x\in\Z,\,t\ge0\}$ whose law is assumed invariant and ergodic under space-time shifts. We prove a quenched invariance principle for the random walk under the minimal moment conditions on the environment; namely, assuming only that the conductances possess the first positive and negative moments. A novel ingredient is the representation of the parabolic coordinates and the corrector via a dual random walk which is considerably easier to analyze. 
\end{quote}

\bigskip
\centerline{\it Dedicated to Jean-Dominique Deuschel}

\section{Introduction}
\nopagebreak
\noindent
The aim of this work is to describe the long-time behavior of a random walk among dynamical random conductances. This problem has enjoyed considerable attention in recent years; we will comment on the relevant literature as soon as the key concepts have been introduced. Throughout this paper we will focus only on one specific instance; namely, the nearest-neighbor random walks on~$\Z$. Our aim is to prove that this walk scales to a non-degenerate Brownian motion assuming only minimal moment conditions on the random environment.

Let us introduce the problem in more precise terms. The aforementioned random ``walk'' is actually a continuous-time Markov chain on~$\Z$ whose dynamics is best described by the (time-dependent) generator~$L_t$ that acts on $f\colon\Z\to\R$ via
\begin{equation}
\label{E:1.1}
(L_t f)(x):=\sum_{z=\pm1}a_t(x,x+z)\bigl[f(x+z)-f(x)\bigr],\quad x\in\Z.
\end{equation}
Here $\{a_t(x,x\pm1)\colon x\in\Z,\,t\ge0\}$ is a family of positive (and finite) numbers that are assumed to obey the symmetry condition
\begin{equation}
a_t(x,x+1)=a_t(x+1,x),\quad x\in\Z,\,\,t\ge0.
\end{equation}
We will refer to $a_t(e)$, for $e=(x,x+1)$, as the \emph{conductance} of edge~$e$ at time~$t$. We will assume that the conductances are defined for all real-valued~$t$ and that they are random, meaning that each $a_t(e)$ is a function of some~$\omega\in\Omega$ in a probability space~$(\Omega,\FF,\BbbP)$. Writing~$\BB(\R)$ for the Borel $\sigma$-algebra on~$\R$, we impose:

\begin{assumption}
\label{ass1}
For each edge~$e$, the map $t,\omega\mapsto a_t(e)$ on~$\R\times\Omega$ is positive, $\BB(\R)\otimes\FF$-measurable, and locally Lebesgue-integrable in~$t$. Moreover, there is a family of space-time shifts, $\tau_{t,x}\colon\Omega\to\Omega$ indexed by $t\in\R$ and~$x\in\Z$, such that
\begin{equation}
\label{E:1.3a}
a_t(x,x+1)\circ\tau_{s,y}=a_{t+s}(x+y,x+y+1),\quad t,s\in\R,\,x,y\in\Z.
\end{equation}
The law~$\BbbP$ is invariant and ergodic with respect to $\{\tau_{t,x}\colon t\in\R,\,x\in\Z\}$.
\end{assumption}

A natural way to interpret the random-walk dynamics is via a Poisson-clock environment: Given a sample of $\{a_t(x,x+1)\colon x\in\Z,\,t\in\R\}$, each edge $e=(x,x+1)$ is endowed with an independent time-inhomogeneous Poisson point process of intensity measure $a_t(e)\textd t$. The above assumptions ensure that this process exists and that no two arrivals, to be called ``rings,'' occur at the same time. The random-walk path is then a \emph{deterministic} function of the Poisson environment: the walk stays at a vertex until an incident edge receives the next ``ring'' at which point it moves to the corresponding neighbor. See Fig.~\ref{fig0} below.

Implementing the Poisson-clock representation rigorously requires showing that the minimal positive solution to the Kolmogorov Backward Equation is non-explosive; i.e., that the number of steps taken by the walk is finite a.s.\ in any finite time. This follows by the assumed local-integrability, stationarity and the Ergodic Theorem. Indeed, for each~$t>0$ there is a (possibly random)~$M\in(0,\infty)$ and a positive density of edges~$e$ (in both lattice directions) where the total jump rate $\int_0^t a_s(e) \textd s$ is bounded by~$M$. Consequently, there is a positive density of edges that receive no ``ring'' in the time-interval~$[0,t]$. Up to time~$t$, the walk is thus effectively confined to a finite set of vertices where the total number of available clock ``rings'' is finite a.s.\ as well. 

Throughout the rest of the paper, we will use the following notation:
\begin{enumerate}
\item[(1)] $X=\{X_t\colon t\ge0\}$ denotes a sample of the above random walk,
\item[(2)] $P^x_a$ denotes the law of~$X$ in a given configuration $a=\{a_t(x,x\pm1)\colon x\in\Z,\,t\in\R\}$ of the conductances subject to the initial condition $P^x_a(X_0=x)=1$, and
\item[(3)] $\E$ denotes expectation with respect to~$\BbbP$.
\end{enumerate}
Our main result is then:

\begin{theorem}[Quenched invariance principle]
\label{thm-main}
Suppose that, on top of Assumption~\ref{ass1}, the conductance law obeys the moment conditions
\begin{equation}
\label{E:1.3}
\E \bigl[a_0(e)\bigr]<\infty\quad\text{and}\quad \E \bigl[a_0(e)^{-1}\bigr]<\infty
\end{equation}
at some (and thus every) edge~$e$. Then there is a constant $\sigma\in(0,\infty)$ such that for any~$T>0$ and~$\BbbP$-a.e.\ sample~$a=\{a_t(x,x+1)\colon x\in\Z,\,t\in\R\}$ of the conductances, the law of
\begin{equation}
X^{(n)}_t:=\frac1{\sqrt n}\,X_{nt},\quad 0\le t\le T,
\end{equation}
induced by~$P^0_a$ on the Skorohod space $\DD[0,T]$ of c\`adl\`ag paths converges, as~$n\to\infty$, weakly to the law of the Brownian motion $\{B_t\colon t\ge0\}$ with $EB_t=0$ and $E(B_t^2)=\sigma^2 t$.
\end{theorem}

Theorem~\ref{thm-main} improves on earlier work by Deuschel and Slowik~\cite{DS16} where the validity of a quenched invariance principle for the corresponding random walk on~$\Z$ was established under the following moment conditions:
\begin{equation}
\label{E:1.6}
\exists p,q\in[1,\infty)\colon\quad
\begin{cases}
&\E \bigl[a_0(e)^p\bigr]<\infty\quad\text{and}\quad \E \bigl[a_0(e)^{-q}\bigr]<\infty
\\*[2mm]
&\displaystyle\frac1{p-1}+\frac1{q(p-1)}<1
\end{cases}
\end{equation}
The algebraic restriction on~$p$ and~$q$ in \eqref{E:1.6} arises from the method of proof which invokes elliptic regularity techniques to construct, and prove sublinearity of, the so called \emph{corrector}, a key object underlying many invariance principles proved so far in this setting. The corresponding problem on~$\Z^d$ for $d\ge2$ has been treated in Andres, Chiarini, Deuschel and Slowik~\cite{ACDS16} albeit under a somewhat different functional relation between~$p$ and~$q$ (and~$d$) than \eqref{E:1.6} might suggest (see~\cite[Remark~1.9]{DS16}).

Although our proof is based on corrector techniques as well, we are able to utilize the one-dimen\-sional nature of the walk to work solely under the weaker conditions \eqref{E:1.3} than \eqref{E:1.6}. Our approach is rooted in that for two-dimensional \emph{static} environments, where a  quenched invariance principle is known to hold under \eqref{E:1.3} in $d=1,2$ (Biskup~\cite{B11}) while requiring $1/p+1/q<2/d$ in~$d\ge3$ (Andres, Deuschel and Slowik~\cite{ASD15}). The need for higher moments in higher dimension has a good reason: for every $p,q\ge1$ satisfying $1/p+1/q>2/(d-1)$, a static environment exists satisfying the moment conditions in \eqref{E:1.3} where the sublinearity of the corrector fails (Biskup and Kumagai~\cite{BK-unpublished}). Whether a quenched invariance principle itself holds just under \eqref{E:1.3} in all~$d\ge1$ remains a subject of extensive debate among experts.

\nopagebreak
\begin{figure}[t]
\vglue1mm
\centerline{\includegraphics[width=0.95\textwidth]{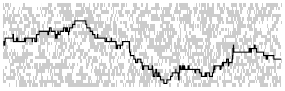}
}
\vglue0mm
\begin{quote}
\small 
\caption{
\label{fig0}
\small
A path of the random walk~$X$ in the dynamical environment composed of independent ON/OFF processes (i.e., the conductances taking values 1 and 0, respectively) with Poisson arrivals of constant (and equal) intensity. The time axis runs horizontally; the shaded areas mark the space-time positions where the corresponding edge is ON. The walk can effectively make a step only across the edges that are ON at that moment of time.}
\normalsize
\end{quote}
\end{figure}

\section{Remarks and outline}
\noindent
We proceed with a couple of remarks. First, the reader may wonder whether the conditions \eqref{E:1.3} are in fact necessary for the result to hold. This is certainly not true for static environments where, thanks to an explicit form of the corrector (see, e.g., Biskup and Prescott~\cite[Introduction]{BP07}) and the fact that we deal with the variable speed random walk (see Barlow and Deuschel~\cite[Theorem~1.1]{BD12} for changes in the constant-speed case), the first condition in \eqref{E:1.3} can be replaced by~$a_0(e)<\infty$ a.s. In the absence of the second condition in \eqref{E:1.3} we actually get a trivial result:

\begin{theorem}[Role of the lower moment condition]
\label{lemma-lower}
Let~$\BbbP$ be be the law of static conductances $\{a(x,x+1)\colon x\in\Z\}$ that are stationary and ergodic with respect to shifts and obey 
\begin{equation}
\BbbP\bigl(a(0,1)<\infty\bigr)=1\quad\text{and}\quad \E\bigl(a(0,x)^{-1}\bigr)=\infty.
\end{equation}
Then for each~$\delta>0$,
\begin{equation}
\label{E:1.6ai}
\E\,P^0_a\bigl(|X_t|\ge\delta\sqrt t\bigr)\,\,\underset{t\to\infty}{\longrightarrow}\,\,0
\end{equation}
In particular, under the diffusive scaling the random walk tends to a vanishing limiting process, at least in the sense of finite-dimensional distributions averaged over the environment.
\end{theorem}

As should be intuitively clear, the main role of the upper moment condition is to prevent blow ups. Here it suffices to consider spatially-homogeneous (dynamical) random environments:

\begin{theorem}[Role of the upper moment condition]
\label{thm-2.2}
Given a stationary ergodic process $\{\eta_t\colon t\in\R\}$ on~$(0,\infty)$ with law~$\BbbP$, define the dynamical conductances via
\begin{equation}
\label{E:2.3o}
a_t(x,x+1):=\eta_t,\quad x\in\Z.
\end{equation}
If $\E\eta_0=\infty$, then for any~$t>0$ and for~$\BbbP$-a.e.\ sample of the conductances, the random variables $\{n^{-1/2}X_{nt}\colon n\ge0\}$ are not tight under~$P^0_a$.
\end{theorem}

These examples show that our moment conditions \eqref{E:1.3} are not only sufficient, but also necessary for a quenched invariance principle with a non-trivial limit process to hold in all the environments satisfying Assumption~\ref{ass1}.

\smallskip
Our second remark concerns the situation when we actually allow the conductances to vanish over sets of times of positive Lebesgue measure. This has been addressed by Biskup and Rodriguez~\cite{BR18}, albeit only in~$d\ge2$, by requiring sufficiently high (namely, $4d+\epsilon$) moments of the quantity
\begin{equation}
T_e:=\inf\Bigl\{t\ge0\colon\int_0^t\textd s\,\,a_s(e)\ge1\Bigr\}.
\end{equation}
We believe that the arguments presented here can be extended to cover the $d=1$ case as well although it is not clear what the minimal moment conditions on~$T_e$ should be. Note that this setting includes some relevant examples; e.g., the random walk on dynamical bond percolation (see Fig.~\ref{fig0}).

Our third remark concerns the dual random walk, which underlies the proofs in the rest of this paper. Leaving the introduction of this walk to Section~\ref{sec3}, we just note that this walk has the same diffusive constant as the main walk of concern in this paper (see Remark~\ref{rem-8.4} for details). It would be of interest to see if a closer --- ideally, path-wise coupling --- relation between these processes could be established. Related to this is the fact that the current proof relies also quite heavily on the assumption that the jumps are only between the nearest neighbors.

Our final remark concerns the fact that the random walk is of variable speed. Here we note that, unlike the case of static environments, in dynamical environments different ways to assign speed --- i.e., normalize the generator --- cannot be related by a time change of the underlying process. At this point, all the existing studies of invariance principles in these cases (namely, the aforementioned references~\cite{ACDS16,BR18}) are restricted to the variable speed case. It is thus of interest to see whether the present approach can be extended to include other versions, most notably discrete-time, as well.

\smallskip
The remainder of this note is organized as follows. In Section~\ref{sec2} we present the standard homogenization argument that gives the convergence in Theorem~\ref{thm-main} subject to two technical claims: existence and sublinearity of the corrector. The main novel contribution of the paper is explained in Sections~\ref{sec3}--\ref{sec4} where we introduce an auxiliary random walk that drives various computation in the rest of the argument. The proof of the technical claims is relegated to Sections~\ref{sec5}--\ref{sec7}. Theorems~\ref{lemma-lower}--\ref{thm-2.2} are proved in Section~\ref{sec8}.

\section{Homogenization argument}
\label{sec2}\nopagebreak\noindent
We are now ready to start discussing the proof of our main results. The argument for convergence builds on well-known techniques from homogenization theory (see Kumagai~\cite{Kumagai} and Biskup~\cite{B11} for recent overviews) which we will explain next. It is the proof of the key technical ingredients --- namely, the existence and sublinearity of the corrector --- that requires a model-specific, and quite non-standard, approach.  

We will henceforth abbreviate
\begin{equation}
b_t(x):=a_t(x,x+1)
\end{equation}
and note that \eqref{E:1.3a} becomes
\begin{equation}
\label{E:b-shift}
b_s(y)\circ\tau_{t,x} = b_{s+t}(y+x),\quad s,t\in\R,\,x,y\in\Z.
\end{equation}
The first point to note is that the structure of the underlying Markov chain gives us the standard ``point of view of the particle:''

\begin{lemma}[Point of view of the particle]
\label{lemma-POVP}
Suppose Assumption~\ref{ass1} holds. Given a sample $a:=\{b_t(x)\colon t\in\R,\,x\in \Z\}$ from~$\BbbP$, let~$\{X_t\colon t\ge0\}$ be a sample from $P^{0}_a$. Then $t\mapsto\tau_{t,X_t}(a)$ is a Markov process on~$\Omega$ with invariant distribution~$\BbbP$. Moreover, the process is ergodic in the sense that, for any $f\in L^1(\BbbP)$,
\begin{equation}
\frac1T\int_0^T \textd t\,\, f\circ\tau_{t,X_t}\,\underset{T\to\infty}\longrightarrow\,\E f
\end{equation}
for $\BbbP$-a.e.~$a\in\Omega$ and $P^0_a$-a.e.~$\{X_t\colon t\ge0\}$.
\end{lemma}

\begin{proofsect}{Proof}
This is standard; see, e.g., Biskup and Rodriguez~\cite[Lemma~4.8]{BR18}.
\end{proofsect}

\begin{figure}[t]
\vglue1mm
\centerline{\includegraphics[width=0.9\textwidth]{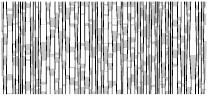}
}
\vglue0mm
\begin{quote}
\small 
\caption{
\label{fig1}
\small
A sample of the parabolic coordinates for the dynamical environment in Fig.~\ref{fig0}. The dark lines mark the ``trajectories'' $t\mapsto\psi(t,x)$ of individual vertices. The time axis runs upwards this time.}
\normalsize
\end{quote}
\end{figure}

Next we introduce the corrector method which relies on the concept of the \emph{parabolic coordinates}. These can be thought of as a time-dependent random embedding of~$\Z$ into~$\R$ that turns the random walk into a martingale; see Fig.~\ref{fig1}. Note that, in static environments, the corresponding object solves a Laplace equation for the generator of the Markov chain and can thus be called a \emph{harmonic coordinate}. In dynamical environments, the Laplace equation is replaced by a parabolic problem; namely, the (reversed-time) heat equation.

\nopagebreak

Recall our notation~$L_t$ for the time dependent generator in \eqref{E:1.1}. The existence and relevant properties of the parabolic coordinates are then the content of: 

\begin{theorem}[Parabolic coordinates]
\label{thm-psi}
Under Assumption~\ref{ass1} and \eqref{E:1.3}, for~$\BbbP$-a.e.\ sample of the conductances there is a map $\psi\colon\R\times\Z\to\R$  such that the following holds:
\settowidth{\leftmargini}{(11)}
\begin{enumerate}
\item[(1)] $t,x\mapsto\psi(t,x)$ is a weak solution to
\begin{equation}
\label{E:2.1}
\frac\partial{\partial t}\psi(t,x)+L_t\psi(t,x)=0,\quad t\in\R,\,x\in\Z,
\end{equation}
with the ``initial'' data $\psi(0,0)=0$.
Moreover, $t\mapsto\psi(t,x)$ is continuous for each~$x\in\Z$.
\item[(2)] For each~$t,s\in\R$ and each~$x,y\in\Z$, the cocycle condition holds
\begin{equation}
\label{E:2.5}
\psi(t+s,x+y)-\psi(t,x)=\psi(s,y)\circ\tau_{t,x}.
\end{equation}
\item[(3)] $\psi(t,x)$ is a jointly measurable function of $t$ and the environment with
\begin{equation}
\label{E:2.6a}
\E\psi(t,x)=x,\quad t\in\R,\,x\in\Z,
\end{equation}
and
\begin{equation}
\label{E:2.7}
\E\bigl(b_0(0)\psi(0,1)^2\bigr)<\infty.
\end{equation}
\end{enumerate}
\item[(4)] Finally, the spatial gradients of~$\psi(t,\cdot)$ are a.s.\ positive,
\begin{equation}
\label{E:2.8}
\psi(t,x+1)-\psi(t,x)>0,\quad t\in\R,\,x\in\Z.
\end{equation}
\end{theorem}

Note that condition \eqref{E:2.8} ensures that under the embedding of~$\Z$ using the parabolic coordinates, the vertices do not swap their order (or, in other words, their space-time trajectories never cross; see Fig.~\ref{fig1}). Deferring the proof to later, we note:

\begin{corollary}
Let $\FF_t:=\sigma(X_s\colon 0\le s\le t)$ and, given a sample~$\{X_t\colon t\ge0\}$, let
\begin{equation}
\label{E:2.9}
M_t:=\psi(t,X_t),\quad t\ge0.
\end{equation}
Then $\{M_t,\FF_t\colon t\ge0\}$ is an $L^2$-martingale with c\`adl\`ag paths and the variance process
\begin{equation}
\langle M\rangle_t:=\int_0^t \textd s\,\,\Theta\circ\tau_{s,X_s}\,,
\end{equation}
where
\begin{equation}
\Theta:=b_0(0)\psi(0,1)^2+b_0(-1)\psi(0,-1)^2.
\end{equation}
\end{corollary}

\begin{proofsect}{Proof}
The continuity of~$t\mapsto \psi(t,x)$ along with the c\`adl\`ag property of~$t\mapsto X_t$ ensure the c\`adl\`ag property of~$t\mapsto M_t$. Recalling that~$X$ has piecewise constant paths a.s., let $\wt N(t)$ denote the number of jumps of~$X$ in the time interval~$[0,t]$. Integrating \eqref{E:2.1} yields
\begin{equation}
\label{E:3.12ui}
M_t=M_0+\int_{(0,t]}\wt N(\textd s)\,\bigl[\psi(s,X_s)-\psi(s,X_{s-})\bigr]-\int_0^t \textd s\,\,(L_s\psi)(s,X_s).
\end{equation}
Since, as~$\epsilon\downarrow0$,
\begin{equation}
E\biggl(\int_t^{t+\epsilon}\wt N(\textd s)\bigl[\psi(s,X_s)-\psi(s,X_{s-})\bigr]\,\bigg|\,\FF_t\biggr)=\int_t^{t+\epsilon}\textd s\,\,(L_s\psi)(s,X_s)+o(\epsilon),
\end{equation}
this shows that $\{M_t,\FF_t\colon t\ge0\}$ is a local martingale. The compensator on the right-hand side of \eqref{E:3.12ui} is (Lebesgue) differentiable, and so the quadratic variation process $[M]$ (using Helland's~\cite{Helland} notation) of~$M$ is carried entirely by its discontinuous part,
\begin{equation}
[M]_t=\int_{(0,t]}\wt N(\textd s)\,\bigl[\psi(s,X_s)-\psi(s,X_{s-})\bigr]^2.
\end{equation}
The variance process~$\langle M\rangle$ is the compensator that makes~$[M]$ a martingale. (We use the cocycle conditions \eqref{E:2.5} to write $\langle M\rangle_t$ using the space-time shifts.)
The condition \eqref{E:2.7} (and the fact that~$\Theta\ge0$) ensures that $t\mapsto\Theta\circ\tau_{t,X_t}$ is locally integrable and, using an elementary localization argument, $M_t$ is thus square integrable for all~$t\ge0$.
\end{proofsect}

As noted before, $x\mapsto\psi(t,x)$ can be thought of as a time-dependent, random embedding of the lattice~$\Z$ into~$\R$ that makes the random walk a martingale. The deformation caused by the change of the embedding, 
\begin{equation}
\chi(t,x):=\psi(t,x)-x
\end{equation}
is the aforementioned \emph{corrector}. A key issue to address now is how much the deformation affects the random walk at the diffusive space-time scales. For this we need:

\begin{theorem}[Sublinearity in diffusive boxes]
\label{thm-chi}
Under Assumption~\ref{ass1} and \eqref{E:1.3}
\begin{equation}
\label{E:2.14}
\max_{\begin{subarray}{c}
x\in\Z\\|x|\le \sqrt n
\end{subarray}}
\,\,
\sup_{\begin{subarray}{c}
t\in\R\\0\le t\le n
\end{subarray}}
\frac{|\chi(t,x)|}{\sqrt n}
\,\underset{n\to\infty}\longrightarrow\,0,\quad\BbbP\text{\rm-a.s.}
\end{equation}
\end{theorem}

The proof of this theorem will be given in Sections~\ref{sec6}-\ref{sec7}. With the help of the above theorems, we can now give:

\begin{proofsect}{Proof of Theorem~\ref{thm-main} from Theorems~\ref{thm-psi}-\ref{thm-chi}}
The following argument is standard; we include it merely for completeness of the exposition.
Consider the martingale~$M$ from \eqref{E:2.9} and let $\langle M\rangle_t$ be its variance process. Lemma~\ref{lemma-POVP} ensures that, for~$\BbbP$-a.e.\ sample of the environment and~$P^0_a$-a.e.\ path of the Markov chain,
\begin{equation}
\label{E:3.17ia}
\frac1t\langle M\rangle_t\,\underset{t\to\infty}\longrightarrow\,\sigma^2:=\E\Theta=2\E\bigl(b_0(0)\psi(0,1)^2\bigr).
\end{equation}
Next recall that~$\wt N(t)$ denotes the number of jumps of~$X$ in time interval~$[0,t]$ and consider the truncated quadratic variation process (using again the notation of Helland~\cite[formula (4.6)]{Helland})
\begin{equation}
\sigma^\epsilon[M](t):=\int_{(0,t]}\wt N(\textd s)\,\bigl[\psi(s,X_s)-\psi(s,X_{s-})\bigr]^2\,1_{\{|\psi(s,X_s)-\psi(s,X_{s-})|>\epsilon\sqrt t\}}.
\end{equation}
By the cocycle conditions \eqref{E:2.5}, the compensator of $\sigma^\epsilon[M]$ is given by
\begin{equation}
\wt\sigma^\epsilon[M](t)=\int_0^t\textd s\,\,\Theta_{\epsilon\sqrt t}\circ\tau_{s,X_s}\,,
\end{equation}
where we set, for general~$r>0$,
\begin{equation}
\Theta_{r}:=b_0(0)\psi(0,1)^2\,1_{\{|\psi(0,1)|>r\}}+b_0(-1)\psi(0,-1)^2\,1_{\{|\psi(0,-1)|>r\}}.
\end{equation}
The Dominated Convergence ensures that $\E\Theta_r\to0$ as~$r\to\infty$. By Lemma~\ref{lemma-POVP} and the downward monotonicity of~$r\mapsto\Theta_r$, for~$\BbbP$-a.e.\ sample~$a$ of the conductances and $P^0_a$-a.e.\ sample of~$X$ we thus get
\begin{equation}
\label{E:3.21ia}
\frac1t\wt\sigma^\epsilon[M](t)\,\underset{t\to\infty}\longrightarrow\,0,\quad \epsilon>0.
\end{equation}
Formulas \eqref{E:3.17ia} and \eqref{E:3.21ia} establish the conditions of the Martingale Functional CLT from Helland~\cite[Theorem~5.1(a)]{Helland} and so, for each~$T>0$, the law of
\begin{equation}
M_t^{(n)}:=\frac1{\sqrt n}\,M_{nt},\quad t\ge0,
\end{equation}
on $\DD[0,T]$ tends weakly, as $n\to\infty$, to that of a Brownian motion $\{B_t\colon t\ge0\}$ with $EB_t=0$ and $E(B_t^2)=\sigma^2 t$. Clearly, $\sigma^2\in(0,\infty)$ by \eqref{E:2.7} and \eqref{E:2.8}.

In order to prove the corresponding statement for the paths of the Markov chain itself, it suffices to show
\begin{equation}
\label{E:2.17}
\sup_{t\le T}\,\bigl|\,X_t^{(n)}-M_t^{(n)}\bigr|\,\underset{n\to\infty}{\overset{P^0_a}\longrightarrow}\,0,\quad\BbbP\text{-a.s.}
\end{equation}
By Theorem~\ref{thm-chi}, for each~$\epsilon>0$ there is a (random) $K$ with~$\BbbP(K<\infty)=1$ such that
\begin{equation}
\bigl|\chi(t,x)\bigr|\le K+\epsilon\bigl(\sqrt t\,+|x|\bigr),\quad t\ge0,\,x\in\Z.
\end{equation}
For~$\epsilon<1$, the triangle inequality converts this to the pointwise estimate
\begin{equation}
|\,X_t^{(n)}-M_t^{(n)}|\le\frac {K}{(1-\epsilon)\sqrt n}+\frac\epsilon{1-\epsilon}\bigl(\sqrt t+|M_t^{(n)}|\bigr).
\end{equation}
Assuming $\epsilon\le1/2$ this gives
\begin{equation}
\sup_{t\le T}\,\bigl|\,X_t^{(n)}-M_t^{(n)}\bigr| \le \frac{2K}{\sqrt n}+2\epsilon \sqrt T+2\epsilon\sup_{t\le T}|M_t^{(n)}|.
\end{equation}
The weak convergence of $M^{(n)}$ to Brownian motion ensures that $\{\sup_{t\le T}|M_t^{(n)}|\colon n\ge1\}$ is tight. Taking $n\to\infty$ followed by~$\epsilon\downarrow0$ then yields \eqref{E:2.17}, as desired.
\end{proofsect}

\section{Dual random walk}
\label{sec3}\nopagebreak\noindent
The proof of our main result has so far been reduced to Theorems~\ref{thm-psi}-\ref{thm-chi} whose proofs constitute the remainder of this paper.   In prior work (namely, \cite{DS16}) these were proved with the help of elliptic-regularity techniques that require the moment conditions \eqref{E:1.6}. As we only wish to assume \eqref{E:1.3}, we will proceed by methods that are tailored to the underlying one-dimensional, and nearest-neighbor, nature of the problem. 

To explain the main idea, let us start with the existence of the parabolic coordinates. Suppose $\psi$ solves \eqref{E:2.1}. Then, as is readily checked,
\begin{equation}
g(t,x):=\psi(t,x+1)-\psi(t,x)
\end{equation}
obeys
\begin{equation}
\label{E:2.6}
-\frac{\partial}{\partial t}g(t,x)=\LL_t^+ g(t,x),\quad t\in\R,\,x\in\Z,
\end{equation}
where
\begin{equation}
\LL_t^+ f(x):=b_t(x+1)f(x+1)+b_t(x-1)f(x-1)-2b_t(x)f(x)\,.
\end{equation}
Our principal observation, and the reason for using the adjoint-operator notation, is that~$\LL_t^+$ is the adjoint in $\ell^2(\Z)$ of
\begin{equation}
\label{E:2.8a}
\LL_t f(x):= b_t(x)\bigl[f(x+1)+f(x-1)-2f(x)\bigr],
\end{equation}
which is the generator of the (variable-speed) \emph{simple symmetric random walk}~$Y$ with jump rate~$2b_t(x)$ at~$x$ at time~$t$. The minus sign on the left-hand side of \eqref{E:2.6} directs us to run this random walk  \emph{backwards} relative to our current labeling of time; \eqref{E:2.6} is then recognized to be the Kolmogorov Forward Equation associated with~$Y$.

Next we recall the requirement that the gradients of~$\psi$ be stationary with respect to the space-time shifts. Hence we expect that
\begin{equation}
g(t,x)=\varphi\circ\tau_{t,x},\quad t\in\R,\,x\in\Z,
\end{equation}
for some measurable function~$\varphi$ of the conductances only. Assuming $\varphi\in L^1(\BbbP)$, equation \eqref{E:2.6} is then equivalent to the statement that the measure
\begin{equation}
\label{E:3.6ua}
\Q(\textd a):=\varphi(a)\BbbP(\textd a)
\end{equation}
is stationary for the evolution $t\mapsto \tau_{-t,Y_t}(a)$ of $a=\{b_t(x)\colon t\in\R,\,x\in\Z\}$ on~$\Omega$ induced by the random walk~$Y$. This suggests that we first extract a stationary distribution~$\Q$ of the environments using the usual averaging procedure and, assuming we can show $\Q\ll\BbbP$, \emph{define} $\varphi$ as the Radon-Nikodym derivative $\frac{\textd\Q}{\textd \BbbP}$. 

We note that~$\varphi$, once constructed, has to be non-negative and a simple argument based on stationarity even gives $\varphi>0$ $\BbbP$-a.s. This implies equivalence of~$\Q$ with~$\BbbP$ (which we need to convert a.s.\ statements under~$\Q$ to those under~$\BbbP$) as well as the ``trajectory non-crossing'' condition \eqref{E:2.8}. The fact that~$\E\varphi=1$, which will also be shown as part of the construction, then gives sublinearity of the corrector in the spatial direction.

\smallskip
In order to implement the above strategy, a number of technical hurdles have to be overcome. The first of these is the very existence of the random walk~$Y$ which requires care due to the dependence of the jump-rates on the (possibly highly irregular) field of the conductances. Then comes the construction of the invariant measure~$\Q$, and the Radon-Nikodym term~$\varphi$, which will be performed in Section~\ref{sec4}. The proof of Theorem~\ref{thm-psi} comes in Section~\ref{sec5}.

Let us start with the construction of the dual random walk~$Y$.
Proceeding along the lines standard in the theory of continuous-time Markov chains (see, e.g., Liggett~\cite{Liggett-MC}), we will first define the transition function of~$Y$ as the minimal positive solution to the Kolmogorov Backward Equations and then, while proving non-explosivity, construct the actual chain as well. Throughout we will regard the conductance configuration as fixed and subject only to the explicitly stated (deterministic) requirements.

We start by defining a family of non-negative kernels $\cmss K_n(t,x;s,y)$ indexed by integers $n\ge0$ and depending on reals $-\infty<t\le s<\infty$ and vertices $x,y\in\Z$, inductively via the iteration scheme 
\begin{multline}
\label{E:back-Kolm-n}
\qquad
\cmss K_{n+1}(s,x;t,y) := \texte^{-\int_t^s \textd u\,\,2b_u(x)}\delta_{x,y}
\\
+\int_t^s\textd r\,\texte^{-\int_r^s \textd u\,\,2b_u(x)}\,b_r(x)\,\sum_{z=\pm1}\cmss K_n(r,x+z;t,y),
\qquad
\end{multline} 
where we set $\cmss K_{0}(s,x;t,y):=0$. Notice that, compared to the usual notation for transition kernels, the evolution runs backwards in time.

\begin{lemma}
Suppose that $t\mapsto b_t(x)$ is locally Lebesgue integrable for all~$x\in\Z$. Then for all $t<s$ and all $x,y\in\Z$, the (Lebesgue) integrals on the right-hand side of \eqref{E:back-Kolm-n} converge and $n\mapsto\cmss K_n(s,x;t,y)$ is non-decreasing and taking values in~$[0,1]$. In particular,
\begin{equation}
\label{E:3.7}
\cmss K(s,x;t,y):=\lim_{n\to\infty} \cmss K_n(s,x;t,y)
\end{equation}
exists in~$[0,1]$. Moreover, for any $t$ and~$y$ fixed, $s,x\mapsto \cmss K(s,x;t,y)$ is a non-negative solution to the Kolmogorov Backward Equation
\begin{equation}
\label{E:3.8}
\frac\partial{\partial s}\cmss K(s,x;t,y)=\LL_s\cmss K(s,\cdot;t,y)(x),\quad s<t,\,y\in\Z,
\end{equation}
where~$\LL_s$ acts on the first spatial variable on the right-hand side and the $s$-derivative is in the Lebesgue sense. The kernel $\cmss K$ is sub-stochastics in the sense that, for all $s\ge t$ and all~$x\in\Z$,
\begin{equation}
\label{E:3.9}
\sum_{y\in\Z}\cmss K(s,x;t,y)\le1.
\end{equation}
Finally, $\cmss K$ transforms covariantly under the space-time shifts; namely,
\begin{equation}
\label{E:4.2ua}
\cmss K(s,x;t,y)\circ\tau_{u,z} = \cmss K(s+u,x+z;t+u,y+z)
\end{equation}
holds for all $s\ge t$, all~$u\in\R$ and all $x,y,z\in\Z$.
\end{lemma}

\begin{proofsect}{Proof}
As is readily checked by induction, we have $\cmss K_n\ge0$ with $n\mapsto \cmss K_n$ is non-decreasing and, thanks to the integrability of $t\mapsto b_t(x)$, also $\sum_{y\in\Z}\cmss K_n(s,x;t,y)\le1$. The limit in~\eqref{E:3.7} thus exists and obeys \eqref{E:3.9}. Passing the limit inside the integral in \eqref{E:back-Kolm-n} using the Monotone Convergence Theorem and some elementary differentiation proves that~$\cmss K$ solves the integral version of \eqref{E:3.8}. As is checked by induction from \eqref{E:back-Kolm-n} and \eqref{E:b-shift}, equation \eqref{E:4.2ua} holds for~$\cmss K_n$; the limit \eqref{E:3.7} then extends it to~$\cmss K$ as well.
\end{proofsect}

A standard question arising in the above context is whether equality holds in \eqref{E:3.9}. As usual, this will be resolved by interpreting~$\cmss K_n$ as the transition probability for a Markov chain restricted to make at most~$n$ steps; equality in \eqref{E:3.9} is then equivalent to non-explosivity of this chain in finite time. We need the following ingredients:
\begin{enumerate}
\item[(1)] $Z:=$ the discrete-time simple symmetric random walk on~$\Z$, and
\item[(2)] $N:=$  an independent rate-1 Poisson point process.
\end{enumerate}
Let $P^x$ denote the joint law of these objects such that $P^x(Z_0=x)=1$. Aiming to define the desired Markov chain as a suitable time-change of the constant-speed continuous-time simple random walk $t\mapsto Z_{N(t)}$, we first need to prove:

\begin{lemma}[Non-explosion]
\label{lemma-non-explosion}
Suppose that $t\mapsto b_{-t}(x)$ is Borel-measurable and locally Lebes\-gue integrable on $(0,\infty)$ for all $x\in\Z$ and, in addition, that (as a function of~$t$)
\begin{equation}
\label{E:env-cond}
b_{-t}(x)>0 \quad\text{\rm for Lebesgue a.e.~$t>0$\quad and}\quad\int_0^\infty \textd t\,\,b_{-t}(x)=\infty,\quad x\in\Z.
\end{equation}
Then for all~$x\in\Z$ and $P^x$-a.e.\ realization of the processes~$Z$ and~$N$ as above, there is a unique continuous $A:[0,\infty)\to[0,\infty)$ satisfying
\begin{equation}
\label{E:3.10}
A(t)=\int_0^t\textd s\,\,2b_{-s}(Z_{N(A(s))}),\quad t\ge0.
\end{equation}
Moreover, we have $P^x(A(t)<\infty)=1$ for each~$t\ge0$ and each~$x\in\Z$.
In particular,
\begin{equation}
Y_t:=Z_{N(A(t))}
\end{equation}
is well defined for all~$t\ge0$ $P^x$-a.s. and obeys
\begin{equation}
\label{E:3.12}
\cmss K(0,x;-t,y)=P^x\bigl(Y_t=y\bigr),\quad x,y\in\Z,\,t\ge0.
\end{equation}
\end{lemma}

\begin{proofsect}{Proof}
The starting point is to solve \eqref{E:3.10} for~$A$. We will do this by constructing its inverse, to be denoted by~$W$. 
Let~$\tau_0:=0<\tau_1<\tau_2<\dots$ mark the arrival times of the Poisson process~$N$. On~$[\tau_n,\tau_{n+1})$ we have $N(\cdot)=n$ and so we may define $W$ inductively by setting $W(0)=W(\tau_0):=0$ and
\begin{equation}
\label{E:3.13}
\int_{W(\tau_n)}^{W(t)}\textd s\,\,2b_{-s}(Z_n) = t-\tau_n,\quad t\in [\tau_n,\tau_{n+1}]
\end{equation}
for all~$n\ge0$. Here the second condition in \eqref{E:env-cond} forces that~$W(t)<\infty$ for all~$t\ge0$ while the integrability and positivity of~$t\mapsto b_{-t}(x)$ assumed in \eqref{E:env-cond} ensure that $t\mapsto W(t)$ is uniquely defined, strictly increasing, continuous on~$[0,\infty)$.

Next let $W(\infty):=\sup_{t\ge0}W(t)$ and define the inverse of~$W$ by
\begin{equation}
\label{E:3.14}
A(t):=\sup\bigr\{r\ge0\colon W(r)\le t\bigr\},\quad 0\le t<W(\infty).
\end{equation}
Set~$t_n:=W(\tau_n)$ and note that $N(A(s))=n$ for $A(s)\in[\tau_n,\tau_{n+1})$ which is equivalent to $s\in[t_n,t_{n+1})$. Using this in \eqref{E:3.13} (and invoking the continuity of~$A$) shows
\begin{equation}
\label{E:4.17ieu}
\int_{t_n}^t \textd s\,\,2b_{-s}(Z_{N(A(s))})=A(t)-\tau_n,\quad t\in[t_n,t_{n+1}].
\end{equation}
As~$t_0=0$, this yields \eqref{E:3.10} for all $t<W(\infty)$ by elementary resummation.

From \eqref{E:back-Kolm-n} we now inductively check that, for all~$t\ge0$,
\begin{equation}
\label{E:3.16}
\cmss K_n(0,x;-t,y)=P^x\bigl(Y_t=y,\,N(A(t))<n\bigr),\quad n\ge0.
\end{equation}
To get \eqref{E:3.12} we have to show that~$Y$ is non-explosive meaning $P^x(N(A(t))<\infty)=1$ for each~$t\ge0$. By \eqref{E:3.14} this boils down to proving $W(\infty)=\infty$ $P^x$-a.s. 
Noting that~$Z$ is recurrent, there is $P^x$-a.s.\ an infinite sequence~$n_0=0<n_1<n_2<\dots$ enumerating the times with~$Z_{n_k}=x$. Then $Z_{N(A(s))}=x$ for $s\in[t_{n_k},t_{n_{k}+1})$ and so, by \eqref{E:4.17ieu},
\begin{equation}
\int_0^{W(\infty)}\!\textd s\,\,2b_{-s}(x)\ge\sum_{k\ge0}\int_{t_{n_k}}^{t_{n_{k}+1}}\!\!\textd s\,\,2b_{-s}(Z_{N(A(s))})
=\sum_{k\ge0}(\tau_{n_{k}+1}-\tau_{n_k}).
\end{equation}
The sum on the right diverges $P^x$-a.s.\ because $\{\tau_{n_{k}+1}-\tau_{n_k}\colon k\ge0\}$ are i.i.d.\ exponential(1) independent of~$Z$. The local integrability of $s\mapsto b_{-s}(x)$ then forces~$W(\infty)=\infty$ $P^x$-a.s., as desired.
\end{proofsect}

As a consequence we now readily get:

\begin{corollary}
Under the assumptions of Lemma~\ref{lemma-non-explosion}, for each each~$x\in\Z$ and each~$s\in\R$, $t,y\mapsto K(s,x;-t,y)$ is a strong solution to the Kolmogorov Forward Equation
\begin{equation}
\label{E:3.17}
-\frac\partial{\partial t}\cmss K(s,x;t,y)=\LL_t^+\cmss K(s,x;t,\cdot)(y),
\end{equation}
for all~$t\le s$ and all~$y\in\Z$.
\end{corollary}

\begin{proofsect}{Proof}
By a simple translation of the environment (which preserves the conditions of Lemma~\ref{lemma-non-explosion}) it suffices to prove this for~$s:=0$. In this case we have the representation \eqref{E:3.16}.
Decomposing according to the last step of the walk~$Y$ we then get
\begin{multline}
\label{E:3.21}
\quad
\cmss K_{n+1}(0,x;-t,y)=P^x\bigl(Y_t=y,\,N(A(t))=0\bigr)+P^x\bigl(Y_t=y,\,0<N(A(t))<n+1\bigr)
\\
=\texte^{-\int_{-t}^0 \textd u\,\,2b_u(y)}\delta_{x,y}
+\int_{-t}^0\textd r\,\texte^{-\int_{-t}^r \textd u\,\,2b_u(y)}\,\sum_{z=\pm1}b_r(y+z)\cmss K_n(0,x;r,y+z).
\quad
\end{multline}
Taking~$n\to\infty$ and using \eqref{E:3.7} along with the Monotone Convergence Theorem we get that~$\cmss K$ satisfies the integral from of \eqref{E:3.17}.
\end{proofsect}

Note that the proof also yields
\begin{equation}
\label{E:3.17a}
-\frac\partial{\partial t}\cmss K_{n+1}(s,x;t,y)=\LL_t^+\cmss K_n(s,x;t,\cdot)(y),\quad n\ge0,
\end{equation}
which will come handy later.

\section{Invariant measure for dual random walk}
\label{sec4}\nopagebreak\noindent
Moving along the strategy outlined at the beginning of Section~\ref{sec3}, we will now construct an invariant distribution~$\Q$ for the Markov chain $t\mapsto\tau_{-t,Y_t}(a)$ on random environments and thus prove Theorem~\ref{thm-psi}. Throughout we consider Assumption~\ref{ass1} and the moment conditions \eqref{E:1.3} as granted. We leave it to the reader to check that this ensures the condition \eqref{E:env-cond} for a.e.\ sample of the conductances.

A standard way to extract an invariant distribution is to average the indicator of an event~$A$ over a finite-stretch of the Markov chain path initiated from the \emph{a priori} measure, and then take a weak subsequential limit. For such an averaged measure, Tonelli's Theorem, the shift-invariance of~$\BbbP$ and \eqref{E:4.2ua} yield
\begin{equation}
\label{E:4.1}
\begin{aligned}
\Q_T(A):&=\frac1T\int_0^T\,\textd t\,\,\E\Bigl(E^0_a\bigl(1_A\circ\tau_{-t,Y_t}\bigr)\Bigr)
\\
&=
\frac1T\int_0^T\,\textd t\,\,\E\Bigl(\,\sum_{y\in\Z}1_A\circ\tau_{-t,y}\cmss K(0,0;-t,y)\Bigr)
\\
&=\E\biggl(1_A\frac1T\int_0^T\,\textd t\,\,\sum_{y\in\Z}\cmss K(t,y;0,0)\biggr).
\end{aligned}
\end{equation}
Writing $\wt\varphi_T$ for the expression following~$1_A$ in \eqref{E:4.1} gives $\Q_T(\textd a):=\wt\varphi_T(a)\BbbP(\textd a)$. Assuming we can prove tightness, every subsequential weak limit of measures~$\Q_T$ as~$T\to\infty$ will then be invariant for the induced chain $t\mapsto\tau_{-t,Y_t}(a)$. 

\smallskip
We will use the above derivation only as motivation; for our purposes, it will be more convenient to work with~$T$ averaged over an exponential distribution. We thus define our approximate Radon-Nikodym term by
\begin{equation}
\label{E:4.2}
\varphi_{\epsilon}:=\epsilon\int_0^\infty\textd t\,\texte^{-\epsilon t}\sum_{y\in\Z}\cmss K(t,y;0,0).
\end{equation}
A similar calculation as in \eqref{E:4.1} shows, with the help of \eqref{E:4.2ua}, that
\begin{equation}
\label{E:phi-int}
\E(\varphi_\epsilon)=1,\quad \epsilon>0,
\end{equation}
and, in particular, $\varphi_\epsilon<\infty$ a.s. 
The main technical problem is to control the ``mass'' of~$\varphi_\epsilon$ in the limit as~$\epsilon\downarrow0$. This will be done via:

\begin{proposition}[Weighted $L^2$-estimate]
\label{prop-moments}
For each~$\epsilon>0$, we have
\begin{equation}
\label{E:4.15}
b_0(0)\varphi_\epsilon^2\in L^1(\BbbP)
\end{equation}
and, in fact,
\begin{equation}
\label{E:4.6ua}
\E\bigl(b_0(0)\varphi_\epsilon^2\bigr)\le\E\bigl(b_0(0)\bigr).
\end{equation}
\end{proposition}

Before we embark on a formal proof, let us note that a similar kind of weighted-$L^2$ estimate appears in most corrector-based approaches to the random conductance model. Disregarding various convergence issues, it is a consequence of the following argument: Introduce the quantity
\begin{equation}
\label{E:chi-def-formal}
\chi_{\epsilon}:=\int_0^\infty\textd t\,\, \texte^{-\epsilon t}\Bigl(b_t(0)\varphi_{\epsilon}\circ\tau_{t,0}-b_t(-1)\varphi_{\epsilon}\circ\tau_{t,-1}\Bigr)\,.
\end{equation}
Then $\varphi_\epsilon-1$ is the spatial gradient of~$\chi_\epsilon$,
\begin{equation}
\label{E:chi-phi-apprx}
\chi_\epsilon\circ\tau_{0,1}-\chi_\epsilon=\varphi_\epsilon-1.
\end{equation}
Moreover, writing~$L$ for the operator~$L_t$ lifted to the space of environments,
\begin{equation}
\label{L-op}
Lf:=b_0(0)\bigl[f\circ\tau_{0,1}-f\bigr]+b_0(-1)\bigl[f\circ\tau_{0,-1}-f\bigr],
\end{equation}
and denoting by
\begin{equation}
\label{V-drift}
V:=b_0(0)-b_0(-1)
\end{equation}
the local drift at the space-time origin, $\chi_\epsilon$ satisfies the ``massive'' corrector equation
\begin{equation}
\label{E:5.9iue}
\frac{\partial}{\partial t}\chi_\epsilon\circ\tau_{t,x}=(\epsilon-L)\chi_\epsilon\circ\tau_{t,x}-V\circ\tau_{t,x}\,.
\end{equation}
These two facts give~$\chi_\epsilon$ the meaning of an approximate, stationary corrector.
Multiplying \eqref{E:5.9iue} at~$t=0$ and $x=0$ by~$\chi_\epsilon$, taking expectation, using that
\begin{equation}
L\chi_\epsilon+V = b_t(0)\varphi_{\epsilon}\circ\tau_{t,0}-b_t(-1)\varphi_{\epsilon}\circ\tau_{t,-1}
\end{equation}
along with the fact that $\E\frac{\partial}{\partial t}\chi_\epsilon^2\circ\tau_{t,0}=0$ thanks to stationarity of~$\BbbP$
produces the standard identity
\begin{equation}
\label{E:4.16ab}
\epsilon\E(\chi_{\epsilon}^2)+\E\bigl(b_0(0)\varphi_{\epsilon}^2\bigr)=\E\bigl(b_0(0)\varphi_{\epsilon}\bigr),
\end{equation}
which, being a direct consequence of the PDE \eqref{E:5.9iue}, can be thought of as a statement of elliptic regularity.
From \eqref{E:4.16ab} we get \eqref{E:4.6ua} by dropping the first term on the left and applying the Cauchy-Schwarz inequality on the right-hand side.

Of course, the main issue with this formal calculation is that, at this point, we have no \emph{a priori} information on the integrability of (and even convergence of the integral defining)~$\chi_\epsilon$. We will therefore need to introduce an additional truncation and work with averaging over space and time instead of the random environment.

\smallskip
Recall our notation~$\cmss K_n$ for the kernels defined in \eqref{E:back-Kolm-n}. We start by introducing a truncated version of~$\varphi_\epsilon$ via
\begin{equation}
\label{E:4.2a}
\varphi_{\epsilon,n}:=\epsilon\int_0^\infty\textd t\,\texte^{-\epsilon t}\sum_{y\in\Z}\cmss K_n(t,y;0,0),\quad n\ge0.
\end{equation}
Since~$n\mapsto\cmss K_n$ is (pointwise) non-decreasing and tending to~$\cmss K$, we have 
\begin{equation}
\varphi_{\epsilon,n}\le\varphi_\epsilon\quad\text{and so}\quad\E\varphi_{\epsilon,n}\le1,\quad n\ge0,
\end{equation}
with~$\varphi_{\epsilon,n}\uparrow\varphi_\epsilon$ as~$n\to\infty$ thanks to the Monotone Convergence Theorem.
The key reason for introducing the truncated objects is that they are pointwise bounded: Since the random walk~$Y$ makes only nearest-neighbor jumps and the kernel~$\cmss K_n$ involves only trajectories with at most~$n$ jumps, the sum in \eqref{E:4.2a} is effectively reduced to~$|y|\le n$. From~$\cmss K_n\le1$ we then have
\begin{equation}
\label{E:4.8a}
\varphi_{\epsilon,n}\le 2n+1,\quad n\ge0.
\end{equation}
Next we introduce the truncated version of \eqref{E:chi-def-formal},
\begin{equation}
\label{E:chi-def}
\chi_{\epsilon,n}:=\int_0^\infty\textd t\, \texte^{-\epsilon t}\Bigl(b_t(0)\varphi_{\epsilon,n}\circ\tau_{t,0}-b_t(-1)\varphi_{\epsilon,n}\circ\tau_{t,-1}\Bigr).
\end{equation}
Here the integral converges absolutely since $t\mapsto\varphi_{\epsilon,n}\circ\tau_{t,0}$ is continuous, $t\mapsto b_t(x)$ is locally integrable and \eqref{E:4.8a} thus gives
\begin{equation}
\label{E:4.13a}
|\chi_{\epsilon,n}|\le(2n+1)\int_0^\infty \textd t\,\,\texte^{-\epsilon t}\bigl[b_t(0)+b_t(-1)\bigr].
\end{equation}
By the first condition in \eqref{E:1.3} the integral has finite expectation under~$\BbbP$; Tonelli's Theorem then implies that the integral is finite $\BbbP$-a.s.
We now claim a finite-$n$ version of~\eqref{E:chi-phi-apprx}:

\begin{lemma}
\label{lemma-4.1}
For all~$\epsilon>0$ and all~$n\ge0$,
\begin{equation}
\label{E:chi-phi}
\chi_{\epsilon,n}\circ\tau_{0,1}-\chi_{\epsilon,n}=\varphi_{\epsilon,n+1}-1\,.
\end{equation}
\end{lemma}

\begin{proofsect}{Proof}
The shift-covariance of the~$\cmss K_n$ kernel implies, for any~$t>0$, that
\begin{equation}
\varphi_{\epsilon,n}\circ\tau_{t,0}=\epsilon\texte^{\epsilon t}\int_t^\infty\textd u\,\,\texte^{-\epsilon u}\,\sum_{y\in\Z}\cmss K_{n}(u,y;t,0).
\end{equation}
The Kolmogorov Forward Equation \eqref{E:3.17a} then yields
\begin{equation}
\frac{\partial}{\partial t}\varphi_{\epsilon,n+1}\circ\tau_{t,0}=\epsilon(\varphi_{\epsilon,n+1}\circ\tau_{t,0}-1)-\LL^+\varphi_{\epsilon,n}\circ\tau_{t,0}\,,
\end{equation}
where the derivative on the left is in the Lebesgue sense and~$\LL^+$ is the operator~$\LL_t^+$ lifted to the space of environments;
\begin{equation}
\label{E:LL+}
\LL^+ f:=b_0(1)f\circ\tau_{0,1}+b_0(-1)f\circ\tau_{0,-1}-2b_0(0)f.
\end{equation}
The definition \eqref{E:chi-def} now shows
\begin{equation}
\begin{aligned}
\chi_{\epsilon,n}\circ\tau_{0,1}-\chi_{\epsilon,n}
&=\int_0^\infty\textd t\,\texte^{-\epsilon t}\LL^+\varphi_{\epsilon,n}\circ\tau_{t,0}
\\
&=\int_0^\infty\textd t\,\,\Bigl[\,\texte^{-\epsilon t}\epsilon(\varphi_{\epsilon,n+1}-1)\circ\tau_{t,0}-\texte^{-\epsilon t}\frac{\partial}{\partial t}\varphi_{\epsilon,n+1}\circ\tau_{t,0}\Bigr]
\\
&=-\int_0^\infty\textd t\,\frac{\partial}{\partial t}\bigl[\texte^{-\epsilon t}(\varphi_{\epsilon,n+1}-1)\circ\tau_{t,0}\bigr]=\varphi_{\epsilon,n+1}-1,
\end{aligned}
\end{equation}
where we also used that $t\mapsto (\varphi_{\epsilon,n+1}-1)\circ\tau_{t,0}$ is bounded.
\end{proofsect}

In light of \eqref{E:chi-phi}, for the integrand in \eqref{E:chi-def} we now get
\begin{equation}
\label{E:4.18a}
b_t(0)\varphi_{\epsilon,n+1}\circ\tau_{t,x}-b_t(-1)\varphi_{\epsilon,n+1}\circ\tau_{t,x-1}
=(L\chi_{\epsilon,n}+V)\circ\tau_{t,x}\,,
\end{equation}
where $L$ and~$V$ are as in \eqref{L-op} and \eqref{V-drift}.
Using this we readily check:

\begin{lemma}[Corrector equation]
For each~$x\in\Z$ and each~$n\ge0$, $t\mapsto \chi_{\epsilon,n}\circ\tau_{t,x}$ is continuous and Lebesgue differentiable with 
\begin{equation}
\label{E:4.13}
\frac{\partial}{\partial t}\chi_{\epsilon,n+1}\circ\tau_{t,x} = \epsilon\chi_{\epsilon,n+1}\circ\tau_{t,x} -(L\chi_{\epsilon,n}+V)\circ\tau_{t,x}\,.
\end{equation}
\end{lemma}

\begin{proofsect}{Proof}
It suffices to prove the claim for~$x=0$. Pick~$t\in\R$. Invoking \eqref{E:4.18a} in \eqref{E:chi-def}, an elementary change of variables yields
\begin{equation}
\label{E:4.14}
\chi_{\epsilon,n+1}\circ\tau_{t,0}
=\texte^{\epsilon t}\int_{t}^\infty\textd s\,\,\texte^{-\epsilon s}(L\chi_{\epsilon,n}+V)\circ\tau_{s,0}.
\end{equation}
The claim follows by differentiation with respect to~$t$ at~$t=0$.
\end{proofsect}

The stationarity of~$\BbbP$ and the first condition in \eqref{E:1.3} imply
\begin{equation}
\label{E:5.26iou}
t\mapsto\int_0^{1}\textd s\,b_s(x)\circ\tau_{t,0} \quad\text{grows sublinearly in~$t$, $\BbbP$-a.s.}
\end{equation}
The bound \eqref{E:4.13a} then shows that $t\mapsto\chi_{\epsilon,n}\circ\tau_{t,0}$ has sublinear growth as well. 
Equipped with these observations, we are ready to give:

\begin{proofsect}{Proof of Proposition~\ref{prop-moments}}
The proof runs parallel to the argument leading up to \eqref{E:4.16ab} except that we average of space-time rather than the environment. We continue writing $L\chi_{\epsilon,n}+V$ as it is concise, but the reader should replace this by the left-hand side of \eqref{E:4.18a} whenever convenient. 

The starting point is to multiply \eqref{E:4.13} by~$\chi_{\epsilon,n+1}\circ\tau_{t,0}$ and integrate over~$t\in[0,r]$, for some~$r>0$. Relabeling~$n+1$ for~$n$, this yields
\begin{equation}
\label{E:4.21}
\chi_{\epsilon,n}^2\circ\tau_{r,0}-\chi_{\epsilon,n}^2=2\epsilon\int_0^r\textd t\,\,\chi_{\epsilon,n}^2\circ\tau_{t,0}
-2\int_0^r\textd t\,\,\bigl[\chi_{\epsilon,n} (L\chi_{\epsilon,n-1}+V)\bigr]\circ\tau_{t,0}.
\end{equation}
The integrals are finite $\BbbP$-a.s.\ by to the fact that~$t\mapsto\chi_{\epsilon,n}\circ\tau_{t,0}$ is continuous and~$t\mapsto b_t(x)$ is locally integrable $\BbbP$-a.s. Next we multiply both sides by~$\texte^{-r/R}$ and integrate over~$r\ge0$. The resulting integrals converge absolutely thanks to the $\BbbP$-a.s.\ sublinear growth of~$t\mapsto\chi_{\epsilon,n}\circ\tau_{t,0}$. Neglecting the contribution of the second term on the left of \eqref{E:4.21} and combining that of the first term with the corresponding term on the right-hand side then shows
\begin{multline}
\qquad
(2\epsilon R-1)\int_0^\infty\textd t\,\,\texte^{-t/R}\,\chi_{\epsilon,n}^2\circ\tau_{t,0}
\\
-2R\int_0^\infty\textd t\,\,\texte^{-t/R}\,\bigl[\chi_{\epsilon,n} (L\chi_{\epsilon,n-1}+V)\bigr]\circ\tau_{t,0}\le0.
\qquad
\end{multline}
For~$2R>1/\epsilon$ (to be assumed next) we can drop the first term. Summing the resulting inequality over its translates by~$x\in\{0,\dots,R\}$, the identity \eqref{E:4.18a} along with Lemma~\ref{lemma-4.1} and integration by parts show
\begin{equation}
\label{E:4.25}
\sum_{x=0}^{R-1}\int_0^\infty\textd t\,\,\texte^{-t/R}\,\bigl[ (\varphi_{\epsilon,n+1}-1)b_0(0)\varphi_{\epsilon,n}\bigr]\circ\tau_{t,x}\le f_R\circ\tau_{0,R}-f_R\circ\tau_{0,-1}\,,
\end{equation}
where~$f_R$ is a ``boundary term'' given explicitly by
\begin{equation}
f_R:=\int_0^\infty\textd t\,\,\texte^{-t/R}\,\bigl[\chi_{\epsilon,n}b_0(0)\varphi_{\epsilon,n}\bigr]\circ\tau_{t,0}.
\end{equation}
Since we are aiming to control the right-hand side of \eqref{E:4.25} in~$\BbbP$-probability, it suffices to focus on the ~$R\to\infty$ behavior of~$f_R$ alone. By \eqref{E:4.8a} and \eqref{E:4.13a}, this quantity is bounded in absolute value by~$(2n+1)^2$ times $h_R(0)+h_R(1)$ where
\begin{equation}
\label{E:4.25a}
h_R(x):=\int_0^\infty \textd t\,\,\texte^{-t/R}\,b_t(0)\biggl(\,\int_0^\infty\textd u\,\,\texte^{-\epsilon u}\,b_{u}(x)\biggr)\circ\tau_{t,0}\,.
\end{equation}
In light of \eqref{E:5.26iou}, the part of the integrand in the large parentheses grows sublinearly in~$t$ a.s.\ Plugging that in, bounding $t\texte^{-t/R}$ by a constant times~$R\texte^{-t/(2R)}$ and noting that, by, say, the $L^1$-part of the Pointwise Ergodic Theorem, the integral of $t\mapsto \texte^{-t/(2R)}b_t(0)$ over all~$t\ge0$ as at most order-$R$ in probability, we get that~$h_R$ and thus also \eqref{E:4.25} are~$o(R^2)$ in probability. From $\varphi_{\epsilon,n}\le\varphi_{\epsilon,n+1}$ we then get
\begin{multline}
\label{E:4.26}
\qquad
\sum_{x=0}^{R-1}\int_0^\infty\textd t\,\,\texte^{-t/R}\,\bigl[ b_0(0)\varphi_{\epsilon,n}^2\bigr]\circ\tau_{t,x}
\\
\le o(R^2)+\sum_{x=0}^{R-1}\int_0^\infty\textd t\,\,\texte^{-t/R}\,\bigl[ b_0(0)\varphi_{\epsilon,n}\bigr]\circ\tau_{t,x}.
\qquad
\end{multline}
The Cauchy-Schwarz inequality bounds the square of the second term on the right by the left-hand side times
\begin{equation}
\sum_{x=0}^{R-1}\int_0^\infty\textd t\,\,\texte^{-t/R}\,b_t(x).
\end{equation}
By the Pointwise Spatial Ergodic Theorem (and our assumptions on~$\BbbP$), this quantity is asymptotic to $R^2\E(b_0(0))$ as~$R\to\infty$ and so
\begin{equation}
\sum_{x=0}^{R-1}\int_0^\infty\textd t\,\,\texte^{-t/R}\,\bigl[ b_0(0)\varphi_{\epsilon,n}^2\bigr]\circ\tau_{t,x}\le  R^2\E(b_0(0))+o(R^2)
\end{equation}
with $o(R^2)/R^2\to0$ in probability as~$R\to\infty$. One more use of the Pointwise Spatial Ergodic Theorem on the left-hand side (which, thanks to the Monotone Convergence Theorem, applies to non-negative random variables even without any moment assumptions) then yields
\begin{equation}
\label{E:eq}
\E\bigl[b_0(0)\varphi_{\epsilon,n}^2\bigr]\le\E\bigl[b_0(0)\bigr].
\end{equation}
The claim now follows from $\varphi_{\epsilon,n}\uparrow\varphi_\epsilon$ and the Monotone Convergence Theorem.
\end{proofsect}

\begin{remark}
Once we have \eqref{E:eq}, the Cauchy-Schwarz inequality along with the first condition in \eqref{E:1.3} show that also~$\chi_\epsilon\in L^2(\BbbP)$. The argument leading up to \eqref{E:4.16ab} can then be applied thus proving the identity \eqref{E:4.16ab} directly. 
\end{remark}

With the weighted-$L^2$ estimate in hand, we can move to the construction of the Radon-Nikodym term~$\varphi$. Instead of working with invariant measures, we proceed by (equivalent) functional-analytic arguments. Consider the linear functional
\begin{equation}
\phi_\epsilon(f):=\E\bigl(\varphi_\epsilon f\bigr),\quad f\in L^\infty(\BbbP),
\end{equation}
and note that it is positive and normalized in the sense that
\begin{equation}
\label{E:4.31}
\phi_\epsilon(f)\le\phi_\epsilon(g)\quad\text{if}\quad f\le g\quad\text{and}\quad\phi_\epsilon(1)=1.
\end{equation}
Writing~$L^0(\BbbP)$ for the set of equivalence classes of measurable functions of the environment, the main outcome of the present section is now:

\begin{theorem}
\label{thm-4.4}
For each~$\epsilon>0$, the linear functional|$\phi_\epsilon$ extends to a continuous linear functional on
\begin{equation}
\HH:=\bigl\{f\in L^0(\BbbP)\colon\E(b_0(0)^{-1}f^2)<\infty\bigr\}
\end{equation}
with norm bounded by $[\E(b_0(0))]^{1/2}$ regardless of~$\epsilon>0$. In particular, weak sequential limits of $\phi_\epsilon$ as~$\epsilon\downarrow0$ exist and take the form $f\mapsto\E(\varphi f)$ for some~$\varphi\in L^0(\BbbP)$ satisfying
\begin{equation}
\label{E:4.32}
\varphi\ge0,\quad \E(\varphi)=1\quad\text{\rm and}\quad \E(b_0(0)\varphi^2)\le\E(b_0(0)).
\end{equation}
In addition, for each~$t>0$ we have
\begin{equation}
\label{E:4.34}
\varphi=\sum_{x\in\Z}\varphi\circ\tau_{t,x}\,\cmss K(t,x;0,0),\quad\BbbP\text{\rm-a.s.}
\end{equation}
In particular, $\varphi$ admits a version such that 
\begin{equation}
\label{E:4.35}
\BbbP(\varphi>0)=1
\end{equation}
and that, on a set of full~$\BbbP$-measure,
$t\mapsto\varphi\circ\tau_{t,x}$ is continuous and weakly differentiable with
\begin{equation}
\label{E:4.36}
\frac{\partial}{\partial t}\varphi\circ\tau_{t,x}+\LL^+\varphi\circ\tau_{t,x}=0,\quad t\in\R,\,x\in\Z,
\end{equation}
where~$\LL^+$ is the operator in \eqref{E:LL+}. The measure~$\Q$ defined from~$\varphi$ via \eqref{E:3.6ua} is stationary and ergodic for the induced Markov chain $t\mapsto\tau_{-t,Y_t}(a)$.
\end{theorem}

\begin{proofsect}{Proof}
Let~$\HH^\star$ denote the space of continuous linear functionals on~$\HH$.
Pick~$f\in L^\infty(\BbbP)$ The Cauchy-Schwarz inequality along with \eqref{E:4.6ua} yield
\begin{equation}
\begin{aligned}
\phi_\epsilon(f)=\E\bigl(\varphi_\epsilon f)
&\le\bigl[\E\bigl(b_0(0)^{-1}f^2\bigr)\bigr]^{1/2}\bigl[\E\bigl(b_0(0)\varphi_\epsilon^2\bigr)\bigr]^{1/2}
\\
&\le\bigl[\E\bigl(b_0(0)^{-1}f^2\bigr)\bigr]^{1/2}\bigl[\E\bigl(b_0(0)\bigr)\bigr]^{1/2}\,.
\end{aligned}
\end{equation}
It follows that~$\phi_\epsilon$ extends continuously to~$\HH$ with the norm bounded by~$[\E(b_0(0))]^{1/2}$. As bounded sequences in~$\HH^\star$ are weakly compact, sequential limits of~$\phi_\epsilon$ as~$\epsilon\downarrow0$ exist and, by the Riesz lemma, take the form~$f\mapsto\E(b_0(0)^{-1} h f)$ for some~$h\in\HH$. Writing $\varphi:=b_0(0)^{-1} h$ we get the second inequality in \eqref{E:4.32}; the equality in \eqref{E:4.32} and non-negativity of~$\varphi$ follow from \eqref{E:4.31} and the fact that $1\in\HH$.

Next we observe that, for any~$t>0$, splitting the integral in \eqref{E:4.2} to an integral over~$[0,t)$ and the other over~$[t,\infty)$, the Chapman-Kolmogorov equations for~$\cmss K$ along with \eqref{E:4.2ua} yield
\begin{equation}
\label{E:4.38}
\varphi_\epsilon
=\epsilon\int_0^t\textd s\,\,\texte^{-\epsilon s}\sum_{x\in\Z}\cmss K(x,s;0,0)
+\texte^{-\epsilon t}\sum_{x\in\Z}
\varphi_\epsilon\circ\tau_{t,x}\cmss K(x,t;0,0).
\end{equation}
The calculation \eqref{E:4.1} shows that the $L^1(\BbbP)$-norm of the first term is $1-\texte^{-\epsilon t}$ which tends to zero as~$\epsilon\downarrow0$. Integrating \eqref{E:4.38} against $f\in L^\infty(\BbbP)$, moving the shift away from~$\varphi$, taking~$\epsilon\downarrow0$ along the sequence where~$\phi_\epsilon$ converges and moving the shift back to~$\varphi$ proves \eqref{E:4.34} with the null set possibly depending on~$t$. 

Now define
\begin{equation}
\label{E:4.38ue}
\overline\varphi:=\int_0^\infty\textd t\,\,\texte^{-t}\Bigl(\,\sum_{x\in\Z}\varphi\circ\tau_{t,x}\,\cmss K(t,x;0,0)\Bigr).
\end{equation}
By \eqref{E:4.34} and Tonelli's Theorem,~$\overline\varphi=\varphi$ $\BbbP$-a.s.\ and so~$\overline\varphi$ is a version of~$\varphi$. As is checked with the help of \eqref{E:4.2a} and a change of variables, $t\mapsto\overline\varphi\circ\tau_{t,x}$ continuous in~$t\in\R$ on a set of full~$\BbbP$-measure. Plugging \eqref{E:4.34} for~$\varphi$ on the right-hand side of \eqref{E:4.38ue} and invoking the Chapman-Kolmogorov equations for~$\cmss K$ shows that \eqref{E:4.34} extends to~$\overline\varphi$ and so we can henceforth regard~$\varphi$ to be this version. 
The differential equation \eqref{E:4.36} now follows from the Kolmogorov Forward Equation \eqref{E:3.17} while \eqref{E:4.35} holds because \eqref{E:4.34} implies
\begin{equation}
\{\varphi=0\}\subseteq\bigl\{\varphi\circ\tau_{t,x}=0\colon x\in\Z,\,t\in\R\bigr\}
\end{equation}
due to strict positivity of~$\cmss K(t,\cdot;0,\cdot)$ for~$t>0$ and \eqref{E:4.34} again. The event on the right is invariant under space-time shift and so, by ergodicity of~$\BbbP$, it has probability zero or one. The case of full measure is ruled out by~$\E(\varphi)=1$. 

The invariance of~$\Q$ for the random walk~$t\mapsto\tau_{-t,Y_t}(a)$ on~$\Omega$ is a consequence of \eqref{E:4.34}. To prove ergodicity, we adapt an argument of Andres~\cite[Proposition~2.1]{A14}. Let~$A$ be a measurable set of environments such that for~$\Q$-a.e.\ $a\in A$ and each~$t>0$ we have~$\tau_{-t,Y_t}(a)\in A$ for~$P_a^0$-a.e.\ sample of~$Y$. This implies
\begin{equation}
\begin{aligned}
0=E_{\Q}(1_{A^\cc}1_{A})
&=\E_{\Q} E^0\bigl(1_{A^\cc}1_A\circ\tau_{-t,Y_t}\bigr)
\\
&=\sum_{x\in\Z}\E\Bigl(\varphi1_{A^\cc}1_A\circ\tau_{-t,x}\cmss K(0,0;-t,x)\Bigr)\,.
\end{aligned}
\end{equation}
But $\varphi>0$ and, for~$t>0$, also $\cmss K(0,0;-t,x)>0$ $\BbbP$-a.s.\ and so we get $1_{A^\cc}1_A\circ\tau_{-t,x}=0$ or, equivalently, $1_A\circ\tau_{-t,x}\le 1_A$ $\BbbP$-a.s.\ for each~$t>0$ and each~$x\in\Z$. Swapping the roles of~$A$ and~$A^\cc$ then gives $1_A=1_A\circ\tau_{-t,x}$ $\BbbP$-a.s.\ for each~$t>0$ and each~$x\in\Z$. By shift-ergodicity of~$\BbbP$, we have~$\BbbP(A)\in\{0,1\}$. 
Since~$\Q$ is equivalent to~$\BbbP$, the same applies~to~$\Q(A)$.
\end{proofsect}

\section{Parabolic coordinates}
\label{sec5}\nopagebreak\noindent
Having established the necessary facts pertaining to the dual random walk~$Y$ we now move to the construction of the parabolic coordinates. This proves the first of the two technical theorems underpinning the main convergence result. We then also prepare the ground for proving the second technical claim by developing an alternative representation for the corrector.

Let~$\varphi$ be a quantity constructed in Theorem~\ref{thm-4.4}; we assume that~$\varphi$ is the version that satisfies \eqref{E:4.36} for all~$t$ and~$x$ on a set of full~$\BbbP$-measure. Set
\begin{equation}
\label{E:5.1}
\chi(t,0):=-\int_0^t\textd s\,\bigl(b_s(0)\varphi\circ\tau_{s,0}-b_s(-1)\varphi\circ\tau_{s,-1}\bigr)\,,\quad t\ge0,
\end{equation}
where the integral converges absolutely $\BbbP$-a.s.\ by Tonelli's Theorem and the fact that $b_0(0)\varphi\in L^1(\BbbP)$ as implied by $b_0(0)\varphi^2\in L^1(\BbbP)$ and the second condition in \eqref{E:1.3}. The quantity~$\chi(t,0)$ will serve as the corrector in time~$t$; compare with its precursor in \eqref{E:chi-def-formal}. Remembering that~$\varphi$ should correspond to the spatial gradients of the parabolic coordinate, the cocycle conditions \eqref{E:2.5} dictate that we define
\begin{equation}
\label{E:5.2}
\psi(t,x):=\sum_{k=0}^{x-1}\varphi\circ\tau_{0,k}+\chi(t,0)\circ\tau_{0,x}\,,\quad x\ge0,\,t\ge0.
\end{equation}
The quantities in \eqref{E:5.1}, resp., \eqref{E:5.2} are defined analogously for negative~$t$, resp.,~$x$, by swapping the limits of the integral/sum and changing the overall sign of the expression. With this definition in hand, we are ready to give:

\begin{proofsect}{Proof of Theorem~\ref{thm-psi}}
A similar calculation to that in the proof of Lemma~\ref{lemma-4.1} shows, with the help of the PDE \eqref{E:4.36} obeyed by~$\varphi$, that
\begin{equation}
\chi(t,0)\circ\tau_{0,x+1}-\chi(t,0)\circ\tau_{t,x} = \varphi\circ\tau_{t,x}-\varphi\circ\tau_{0,x}.
\end{equation}
This readily implies
\begin{equation}
\psi(t,x+1)-\psi(t,x)=\varphi\circ\tau_{t,x}
\end{equation}
and proves the cocycle condition \eqref{E:2.5}. The PDE \eqref{E:2.1} obeyed by~$\psi$ is then a direct consequence of the definition \eqref{E:5.1}. The identities \twoeqref{E:2.6a}{E:2.7} follow from \eqref{E:4.32} while \eqref{E:2.8} is a rewrite of~\eqref{E:4.35}.
\end{proofsect}

Although the formula \eqref{E:5.1} serves well for the construction of the parabolic coordinate, it appears less amenable for the purposes of proving Theorem~\ref{thm-chi}. There we will use a different representation which we will prove next:

\begin{proposition}
\label{prop-chi-repr}
For each~$t\ge0$,
\begin{equation}
\label{E:3.1}
\chi(t,0)=
\sum_{x<0}\sum_{y\ge0}\varphi\circ\tau_{t,x}\, \cmss K(t,x;0,y)
-\sum_{x\ge0}\sum_{y<0}\varphi\circ\tau_{t,x}\, \cmss K(t,x;0,y)\,,
\end{equation}
where each of the double sums converges to a finite number $\BbbP$-a.s.
\end{proposition}

For the proof of a.s.\ convergence we first show:

\begin{lemma}
There is a constant $c>0$ such that each~$t\ge0$,
\begin{equation}
\label{E:3.2}
\E\Bigl(\varphi\sum_{z\in\Z}\cmss K(0,0;-t,z)|z|\Bigr)\le c\sqrt t.
\end{equation}
\end{lemma}

\begin{proofsect}{Proof}
Using the stationary distribution on environments, $\Q(\textd a):=\varphi(a)\BbbP(\textd a)$, the quantity in question is recognized as the left-hand side of
\begin{equation}
E_\Q E^0(|Y_t|)\le\bigl[E_\Q E^0(Y_t^2)\bigr]^{1/2}.
\end{equation}
Since~$t\mapsto Y_t$ is a martingale with associated variance process
\begin{equation}
\langle Y\rangle_t = \int_0^\infty \textd s\,\,2b_{-s}(Y_s) = \int_0^\infty \textd s\,\,2b_0(0)\circ\tau_{-s,Y_s}\,,
\end{equation}
from stationarity of~$\Q$ under $t\mapsto\tau_{-t,Y_t}(a)$ we readily get
\begin{equation}
E_\Q E^0(Y_t^2)=E_\Q E^0(\langle Y\rangle_t) = 2t E_\Q\bigl(b_0(0)\bigr)=2t \E\bigl(b_0(0)\varphi\bigr).
\end{equation}
As noted above, the last expectation is finite by \eqref{E:4.6ua}.
\end{proofsect}

\begin{proofsect}{Proof of Proposition~\ref{prop-chi-repr}}
Fix~$t\ge0$. A shift of the environment and a change of variables show that the expectation under~$\BbbP$ of the \emph{sum} of the two terms in \eqref{E:3.1} equals the expectation in \eqref{E:3.2}. Since~$\varphi>0$ $\BbbP$-a.s., the sums in \eqref{E:3.1} converge to a finite number $\BbbP$-a.s. Denoting, with some abuse of our earlier notation, by~$\chi_n(t)$ the quantity in \eqref{E:3.2} with the sums over~$x$ and~$y$ additionally restricted to values in~$[-n,n]$, we in particular have $\chi_n(t)\to\chi(t,0)$ as~$n\to\infty$ a.s.\ by the Dominated Convergence Theorem. 

We will now calculate  the $t$-derivative of~$\chi_n(t)$. Fix~$y\in\Z$ and let us temporarily abbreviate
\begin{equation}
\varphi_x:=\varphi\circ\tau_{t,x},\quad b_x:=b_t(x)\quad\text{and}\quad K_x:=\cmss K(t,x;0,y).
\end{equation}
Then \eqref{E:4.36} reads as
\begin{equation}
\frac{\partial}{\partial t}\varphi\circ\tau_{t,x} = -b_{x+1}\varphi_{x+1}-b_{x-1}\varphi_{x-1}+2b_x\varphi_x\end{equation}
while the Backward Kolmogorov Equation \eqref{E:3.8} becomes
\begin{equation}
\frac{\partial}{\partial t}\cmss K(t,x;0,y)=b_x\bigl[K_{x+1}+K_{x-1}-2K_x\bigr].
\end{equation}
The product rule for the derivative then shows
\begin{multline}
\quad\frac{\partial}{\partial t}\bigl(\varphi\circ\tau_{t,x}\cmss K(t,x;0,y)\bigr)
\\
=\bigl(b_x\varphi_x K_{x+1}-b_{x-1}\varphi_{x-1}K_x\bigr)+\bigl(b_x\varphi_x K_{x-1}-b_{x+1}\varphi_{x+1}K_x\bigr).
\quad
\end{multline}
Using the standard telescoping argument, we have
\begin{equation}
\label{E:6.14o}
\begin{aligned}
\sum_{x=-n}^{-1}\bigl(b_x\varphi_x K_{x+1}-b_{x-1}\varphi_{x-1}K_x\bigr) 
&= b_{-1}\varphi_{-1}K_0-b_{-n-1}\varphi_{-n-1}K_{-n},
\\
\sum_{x=-n}^{-1}
\bigl(b_x\varphi_x K_{x-1}-b_{x+1}\varphi_{x+1}K_x\bigr)
&=b_{-n}\varphi_{-n}K_{-n-1}-b_0\varphi_0 K_{-1}.
\end{aligned}
\end{equation}
Similarly we obtain
\begin{equation}
\label{E:6.15o}
\begin{aligned}
\sum_{x=0}^n\bigl(b_x\varphi_x K_{x+1}-b_{x-1}\varphi_{x-1}K_x\bigr) 
&= b_n\varphi_nK_{n+1}-b_{-1}\varphi_{-1}K_0,
\\
\sum_{x=0}^n
\bigl(b_x\varphi_x K_{x-1}-b_{x+1}\varphi_{x+1}K_x\bigr)
&=b_0\varphi_0 K_{-1}-b_{n+1}\varphi_{n+1}K_n.
\end{aligned}
\end{equation}
We will now return to the full notation while still abbreviating $(b\varphi)_{t,x}:=b_t(x)\varphi\circ\tau_{t,x}$. Summing \twoeqref{E:6.14o}{E:6.15o} over~$y$ in the respective range of values (still confined to~$[-n,n]$) and then subtracting the sums in \eqref{E:6.15o} from those in \eqref{E:6.14o} yields
\begin{equation}
\label{E:6.16o}
\begin{aligned}
\frac{\partial}{\partial t}\chi_n&(t)
=(b\varphi)_{t,-1}\,P^{0}\bigl(|Y_t|\le n\bigr)\circ\tau_{t,0}-(b\varphi)_{t,0}\,P^{-1}\bigl(|Y_t|\le n\bigr)\circ\tau_{t,0}
\\
&-(b\varphi)_{t,-n-1}\,P^{-n}(0\le Y_t\le n)\circ\tau_{t,0}+(b\varphi)_{-n}\,P^{-n-1}(0\le Y_t\le n)\circ\tau_{t,0}
\\
&-(b\varphi)_{t,n}\,P^{n+1}(-n\le Y_t<0)\circ\tau_{t,0}+(b\varphi)_{t,n+1}\,P^{n}(-n\le Y_t<0)\circ\tau_{t,0}.
\end{aligned}
\end{equation}
Here the first term on the right of \eqref{E:6.16o} arose by combining the contributions from the terms $b_0\varphi_n K_{-1}$ in \twoeqref{E:6.14o}{E:6.15o}. Similarly, the second term combines the contributions from the term $b_{-1}\varphi_{-1}K_{0}$. The remaining terms in \eqref{E:6.16o} collect the contributions of the terms $b_{-n-1}\varphi_{-n-1}K_{-n}$, $b_{-n}\varphi_{-n}K_{-n-1}$, $b_n\varphi_n K_{n+1}$ and $b_{n+1}\varphi_{n+1}K_n$, respectively.

The first two terms on the right of \eqref{E:6.16o} dominate the expression in the limit~$n\to\infty$. Indeed, integrating over a compact interval of~$t$ and taking expectation with respect to~$\BbbP$, the remaining four terms on the right of \eqref{E:6.16o} converge to zero in~$L^1(\BbbP)$ as~$n\to\infty$. The term $P^x(|Y_t|\le n)$ in turn increases to one as~$n\to\infty$ for both~$x=0,-1$. The Monotone Convergence Theorem gives
\begin{equation}
\text{r.h.s.\ of \eqref{E:3.1}}=\int_0^t\textd s\,\bigl(b_s(-1)\varphi\circ\tau_{s,-1}-b_s(0)\varphi\circ\tau_{s,0}\bigr),\quad \BbbP\text{-a.s.}
\end{equation}
The quantity on the right is~$\chi(t,0)$, as desired.
\end{proofsect}

\begin{remark}
The reader may wonder at this point how we arrived at the above alternative expression for~$\chi(t,0)$ in the first place. This was done as follows. We know that the spatial gradients of the corrector are given by~$\varphi-1$. Setting
\begin{equation}
\wt\chi_\epsilon:=-\sum_{x\ge0}\frac{\varphi\circ\tau_{0,x}-1}{(1+\epsilon)^{x+1}},
\end{equation}
where the sum converges because $x\mapsto \varphi\circ\tau_{0,x}$ has a sublinear growth,
we then get
\begin{equation}
\wt \chi_\epsilon\circ\tau_{0,1}-\wt\chi_\epsilon=\varphi-1-\epsilon\wt\chi_\epsilon.
\end{equation}
This indicates that $\wt\chi_\epsilon\circ\tau_{t,x}$ is an approximate (stationary) corrector at space time position~$(t,x)$. We should thus be able to approximate $\chi(t,0)$ by the quantity
\begin{equation}
\wt\chi_\epsilon\circ\tau_{t,0}-\wt\chi_\epsilon = 
\sum_{x\ge0}\frac{\varphi\circ\tau_{0,x}-\varphi\circ\tau_{t,x}}{(1+\epsilon)^{x+1}}\,.
\end{equation}
Using \eqref{E:4.34} for the term $\varphi\circ\tau_{0,x}$ and invoking \eqref{E:4.2ua} along with the fact that~$y\mapsto\cmss K(t,x;0,y)$ is a probability mass function, this is recast as 
\begin{equation}
\label{E:6.21o}
\wt\chi_\epsilon\circ\tau_{0,t}-\wt\chi_\epsilon = 
\sum_{x,y\in\Z}\varphi\circ\tau_{t,x}
\biggl[\frac{\cmss K(t,x;0,y)}{(1+\epsilon)^{y+1}}1_{\{y\ge0\}}-
\frac{\cmss K(t,x;0,y)}{(1+\epsilon)^{x+1}}1_{\{x\ge0\}}\biggr]\,.
\end{equation}
Noting that
\begin{equation}
1_{\{y\ge0\}}-1_{\{x\ge0\}}=1_{\{x<0\}}1_{\{y\ge0\}}-1_{\{x\ge0\}}1_{\{y<0\}},
\end{equation}
taking, at least formally, the limit~$\epsilon\downarrow0$ in \eqref{E:6.21o} we then discover \eqref{E:3.1}.
\end{remark}

\section{Corrector sublinearity}
\label{sec6}\nopagebreak\noindent
To make the proof of our main result complete, it remains to establish the everywhere sublinearity of the corrector as stated in Theorem~\ref{thm-chi}. We will proceed by the argument developed in Berger and Biskup~\cite{BB07} for the random walk on two-dimensional supercritical bond-percolation clusters which was later extended (Biskup~\cite{B11}) to random walks in general ergodic conductance models subject to moment conditions of the type \eqref{E:1.3}. A key novel ingredient, stated in Proposition~\ref{prop-6.2}, is proved in Section~\ref{sec7}.

The starting point is sublinearity in the spatial direction:

\begin{lemma}[Sublinearity in space]
\label{lemma-6.1}
For $\BbbP$-a.e.\ sample of the environment,
\begin{equation}
\label{E:6.4}
\lim_{n\to\pm\infty}\frac{|\chi(0,n)|}{n}=0,\quad\BbbP\text{\rm-a.s.}
\end{equation}
\end{lemma}

\begin{proofsect}{Proof}
We follow the proof of \cite[Lemma~4.8]{B11}.
Fix~$n\in\N$ and~$t\ge0$. The cocycle conditions give
\begin{equation}
\chi(0,n)=\sum_{k=0}^{n-1}\chi(0,1)\circ\tau_{0,k}
\end{equation}
and
\begin{equation}
\label{E:6.3}
\chi(0,n)\circ\tau_{t,0}=\chi(0,n)+\chi(0,t)\circ\tau_{0,n}-\chi(0,t).
\end{equation}
Since $\chi(0,1)=\varphi-1$, the equality in \eqref{E:4.32} shows
\begin{equation}
\chi(0,1)\in L^1(\BbbP)\quad\text{and}\quad\E\chi(0,1)=0.
\end{equation}
Birkhoff's Pointwise Ergodic Theorem then gives 
\begin{equation}
\label{E:6.5}
\overline\chi:=\lim_{n\to\infty}\frac{\chi(0,n)}n\quad\text{exists $\BbbP$-a.s.\ and in~$L^1(\BbbP)$}.
\end{equation}
The representation \eqref{E:5.1} along with $b_0(0)\varphi\in L^1(\BbbP)$ ensure
\begin{equation}
\chi(t,0)\in L^1(\BbbP).
\end{equation}
From \eqref{E:6.3} (and $L^1$-convergence) we then get $\overline\chi\circ\tau_{t,0}=\overline\chi$ $\BbbP$-a.s. for each~$t\ge0$. The limit definition in \eqref{E:6.5} ensures that, also, $\overline\chi\circ\tau_{0,x}=\overline\chi$ $\BbbP$-a.s.\ for each~$x\in\Z$. Hence,~$\overline\chi$ is shift invariant, and thus constant~$\BbbP$-a.s.\ by the assumed ergodicity of~$\BbbP$. Using the $L^1$-convergence part of \eqref{E:6.5} we get
\begin{equation}
\overline\chi=\E\overline\chi=\E\chi(0,1)=0.
\end{equation}
This proves the claim for~$n\to\infty$ limit; replacing~$\chi(0,1)$ by~$\chi(0,-1)$ extends this to $n\to-\infty$ limit as well.
\end{proofsect}

Looking at how the ranges of $x$ and~$t$ in \eqref{E:2.14} scale with~$n$, for the behavior of the corrector in time we need to actually prove a subdiffusive growth estimate:

\begin{proposition}[Sudiffusivity in time]
\label{prop-6.2}
For $\BbbP$-a.e.\ sample of the environment,
\begin{equation}
\label{E:6.1}
\lim_{t\to\infty}\frac{|\chi(t,0)|}{\sqrt t}=0.
\end{equation}
\end{proposition}

We remark that finding a representation of the corrector that makes subdiffusivity of the corrector in time transparent has been the primary driving force behind the approach developed in the present paper. Before we delve into its proof (which is deferred to the next section), let us show how it implies the desired theorem:

\begin{proofsect}{Proof of Theorem~\ref{thm-chi} from Proposition~\ref{prop-6.2}}
We follow arguments developed in Berger and Biskup~\cite{BB07}; see also~\cite[Lemma~4.12]{B11}. First we identify a ``good grid'' of space-time points where the corrector can be controlled by way of ergodic-theoretical and geometric arguments. The oscillation of the corrector over the ``holes'' left out by the grid is then controlled by methods of harmonic analysis. The proof is divided into three steps.

\smallskip\noindent
\textsl{Step 1 (Definition of good grid):}
Let us call the space-time point $(0,0)$ $K,\epsilon$-good if
\begin{equation}
\bigl|\chi(0,x)\bigr|\le K+\epsilon|x|,\quad x\in\Z,
\end{equation}
and
\begin{equation}
\bigl|\chi(t,0)\bigr|\le K+\epsilon\sqrt{t},\quad t\ge0.
\end{equation}
Similarly, we will call $(x,t)$ $K,\epsilon$-good in environment~$a$ if $(0,0)$ is $K,\epsilon$-good in the environment $\tau_{t,x}(a)$. In light of \eqref{E:6.4} and \eqref{E:6.1},
\begin{equation}
\label{E:6.11}
\BbbP\bigl((t,x)\text{ is $K,\epsilon$-good}\bigr)\,\underset{K\to\infty}\longrightarrow\,1
\end{equation}
holds for all~$\epsilon>0$ and all~$(t,x)\in[0,\infty)\times\Z$.

Let $\rho_{K,\epsilon}$ be the density of $K,\epsilon$-good points in~$\Z$; this quantity exist by Birkhoff's Ergodic Theorem and is generally random but, since its expectation equals the probability in \eqref{E:6.11}, from the obvious monotonicity in~$K$ we have
\begin{equation}
\rho_{K,\epsilon}\,\underset{K\to\infty}\longrightarrow\,1,\qquad\BbbP\text{-a.s.}
\end{equation}
Similarly, if~$\theta_{K,\epsilon}$ is the density of~$\{n\in\N\colon (n,0)\text{ is $K,\epsilon$-good}\}$ in~$\N$ (remember that~$t\mapsto\chi(t,0)$ is continuous so checking only integer times will be enough) we have\begin{equation}
\theta_{K,\epsilon}\,\underset{K\to\infty}\longrightarrow\,1,\qquad\BbbP\text{-a.s.}
\end{equation}
It follows that, for each~$\epsilon>0$ and~$\BbbP$-a.e.\ environment~$a$ there is $K=K(a)<\infty$ such that
\begin{equation}
\rho_{K,\epsilon}\ge\frac12,\quad \theta_{K,\epsilon}\ge\frac12\quad\text{and}\quad (0,0)\text{ is $K,\epsilon$-good}.
\end{equation}
We now fix this~$K$ and let $\G_{K,\epsilon}$ denote the set of $(t,x)\in[0,\infty)\times\Z$ such that \emph{at least one} of the following conditions holds:
\begin{enumerate}
\item[(1)] $t=0$ or~$x=0$ or both,
\item[(2)] $t$ is integer and~$(t,0)$ is $K,\epsilon$-good,
\item[(3)] $(0,x)$ is $K,\epsilon$-good.
\end{enumerate}
The set~$\G_{K,\epsilon}$ is the aforementioned ``good grid.''

\smallskip\noindent
\textsl{Step 2 (Estimating $\chi$ on good grid):}
We now derive a pointwise estimate of the corrector on the good grid. Note that each~$(t,x)\in\G_{K,\epsilon}$ can be connected to the origin by following a pair of horizontal and vertical lines that lie entirely in~$\G_{K,\epsilon}$ --- which line comes first depends on which of the three condition above applies at~$(t,x)$; one or both lines are trivial when~(1) is in force. Since these line segments meet at a $K,\epsilon$-good point, the cocycle condition and the triangle inequality show
\begin{equation}
\label{E:6.15}
\bigl|\chi(t,x)\bigr|\le 2K+\epsilon|x|+\epsilon\sqrt t,\quad (x,t)\in\G_{K,\epsilon}.
\end{equation}
It remains to control the corrector at points away from~$\G_{K,\epsilon}$. 

Note that the ``holes'' left out by~$\G_{K,\epsilon}$ are rectangles bounded by horizontal and vertical lines in~$\G_{K,\epsilon}$. We will write~$\partial R$ for the points in~$\G_{K,\epsilon}$ bounding rectangle~$R$ (which we think as disjoint from~$\G_{K,\epsilon}$). Next recall that the parabolic coordinates are defined so that~$\psi(t,X_t)$ is a martingale. Using the Optional Stopping Theorem (or the PDE for~$\psi$ directly), this implies a Maximum Principle: For any rectangle~$R$ as above and any~$x_0\in R\cup\partial R$,
\begin{equation}
\sup_{(t,x)\in R}\bigl|\psi(t,x)-x_0\bigr|\le\sup_{(t,x)\in\partial R}\bigl|\psi(t,x)-x_0\bigr|\,.
\end{equation}
Since $\chi(t,x)=\psi(t,x)-x$, we thus get
\begin{equation}
\label{E:6.17}
\sup_{(t,x)\in R}\bigl|\chi(t,x)\bigr|\le\sup_{(t,x)\in\partial R}\bigl|\chi(t,x)\bigr|+\diam_{\Z}(R),
\end{equation}
where~$\diam_{\Z}(R)$ is the diameter of the projection of~$R\cup\partial R$ onto the spatial coordinate. Since~$\partial R\subset\G_{K,\epsilon}$, the supremum on the right can be controlled via \eqref{E:6.15} provided we can control the diameter of any rectangle that intersects $[0,n]\times[-\sqrt n,\sqrt n]$; this then takes care also of the second term on the right. 

\smallskip\noindent
\textsl{Step 3 (Away from good grid):}
Let~$\{x_k\colon k\in\Z\}$, with~$x_0:=0$, be the increasing sequence enumerating $K,\epsilon$-good points on the line~$t=0$; this sequence exists by the fact that~$\rho_{K,\epsilon}>0$ (note that the left and right densities of $\epsilon,K$-good points are equal~$\BbbP$-a.s.). The existence and positivity of the density of good points implies
\begin{equation}
\lim_{k\to\pm\infty}\frac{|x_k-x_{k-1}|}{k}=0.
\end{equation}
Similarly, letting $\{t_k\colon k\ge0\}$, where~$t_0:=0$, enumerate the $K,\epsilon$-good points with integer time coordinate and zero space coordinate, we have
\begin{equation}
\lim_{k\to\infty}\frac{|t_k-t_{k-1}|}k=0.
\end{equation}
Since~$x_k/k$ as well as~$t_k/k$ tend to positive numbers as~$k\to\infty$, there is a (random) $\wt K<\infty$ such that, for all~$k$,
\begin{equation}
|x_k-x_{k-1}|\le\wt K+\epsilon(|x_k|\wedge|x_{k-1}|)\quad\text{and}\quad|t_n-t_{k-1}|\le\wt K+\epsilon t_{k-1}.
\end{equation}
It follows that, once $n\ge\wt K+\epsilon n$, any rectangle~$R$ in~$([0,\infty)\times\Z)\smallsetminus\G_{K,\epsilon}$ that intersects $[0,n]\times[-\sqrt n,\sqrt n]$ satisfies $R\cup\partial R\subset[0,2n]\times[-\sqrt{2n},\sqrt{2n}]$ and
\begin{equation}
\diam_\Z(R)\le\sqrt{\wt K+\epsilon n}\,.
\end{equation}
Combining this with \eqref{E:6.15} and \eqref{E:6.17} yields
\begin{equation}
\sup_{0\le t\le n}\,\,\,\max_{|x|\le \sqrt n}\,\,\,\bigl|\chi(t,x)\bigr|\le 2K+\epsilon \sqrt{2n}+\epsilon\sqrt{2 n}+\sqrt{\wt K+\epsilon n}.
\end{equation}
Dividing by~$\sqrt n$ and taking~$n\to\infty$ followed by~$\epsilon\downarrow0$ then yields the claim.
\end{proofsect}

\section{Subdiffusivity in time}
\label{sec7}\nopagebreak\noindent
As a final point of the proof, it remains to prove the subdiffusive estimate for the corrector in time. It is here where we will benefit from the representation in Proposition~\ref{prop-chi-repr}. As it turns out, it suffices to focus on the limit of large \emph{negative} times. The cocycle conditions give $\chi(-t,0)=-\chi(t,0)\circ\tau_{-t,0}$ for any~$t>0$, and so
\begin{equation}
\label{E:4.2e}
\chi(-t,0)=\sum_{x\ge0}\varphi\circ\tau_{x}\, P^{x}(Y_{t}<0)-\sum_{x<0}\varphi\circ\tau_{0,x}\, P^{x}(Y_{t}\ge 0),\quad t\ge0,
\end{equation}
where~$Y$ is the dual random walk. We start by showing that the sums in \eqref{E:4.2e} are dominated by~$x$-values of order~$\sqrt t$:

\begin{lemma}
\label{lemma-7.1}
For~$\BbbP$-a.e.\ environment, 
\begin{equation}
\label{E:7.8}
\lim_{M\to\infty}\,\limsup_{t\to\infty}\,\frac1{\sqrt t}\sum_{x\ge M\sqrt t}\varphi\circ\tau_{0,x}\,P^x\bigl(Y_t<0\bigr)=0
\end{equation}
and
\begin{equation}
\lim_{M\to\infty}\,\limsup_{t\to\infty}\,\frac1{\sqrt t}\sum_{x\le -M\sqrt t}\varphi\circ\tau_{0,x}\,P^x\bigl(Y_t>0\bigr)=0.
\end{equation}
\end{lemma}

\begin{proofsect}{Proof}
By symmetry it suffices to prove just \eqref{E:7.8}. Recall the definition \eqref{E:3.10} of the time-change process~$A(t)$ that links~$Y$ to the discrete time simple symmetric random walk~$(Z_n)_{n\ge0}$ and an independent rate-1 Poisson process $(N(t))_{t\ge0}$. Pick $p\in(0,1/2)$ and note that the sum in \eqref{E:7.8} is bounded by the sum of the following terms
\begin{equation}
\text{\rm I}_M(t):=\sum_{x\ge M\sqrt t}\varphi\circ\tau_{0,x}\,P^x\Bigl(A(t)\ge x ^{2(1-p)}t^p\Bigr)
\end{equation}
and
\begin{equation}
\label{E:7.14}
\text{\rm II}_M(t):=\sum_{x\ge M\sqrt t}\varphi\circ\tau_{0,x}\,P^x\Bigl(A(t)\le x^{2(1-p)}t^p,\,Y_t<0\Bigr).
\end{equation}
We will now estimate these two terms separately.

Since $\textd A(t)=2b_0(0)\circ\tau_{-t,Y_t}\textd t$, we can analyze the behavior $t\mapsto A(t)$ by following the evolution of the environment from the point of view of the random walk~$Y$. To this end, define the maximal function $A^\star:=\sup_{t>0}\frac{A(t)}t$. The Markov inequality shows, for any~$q>0$, that
\begin{equation}
\begin{aligned}
\text{\rm I}_M(t)
&\le \sum_{x\ge M\sqrt t}\varphi\circ\tau_{0,x}\,P^x\Bigl(A^\star\ge (x^2/t)^{(1-p)}\Bigr)
\\
&\le t^q\sum_{x\ge M\sqrt t}\varphi\circ\tau_{0,x}\,\frac1{x^{2q}}E^x\bigl((A^\star)^{\frac q{1-p}}\bigr).
\end{aligned}
\end{equation}
Since~$b_0(0)\in L^1(\Q)$ and~$\Q$ is invariant for the environment observed from the walk~$Y$, the Maximal Ergodic Theorem gives~$E_\Q E^0((A^\star)^r)<\infty$ for all~$r\in(0,1)$ and so
\begin{equation}
\label{E:8.7uia}
\forall q\in(0,1-p)\colon \quad\varphi\, E^0\bigl((A^\star)^{\frac q{1-p}}\bigr)\in L^1(\BbbP).
\end{equation}
Now use the fact that if~$f\in L^1(\BbbP)$ is non-negative, and $f^\star:=\sup_{n\ge1}\frac1n\sum_{k=0}^{n-1} f\circ\tau_{0,k}$ is the associated maximal function under spatial shifts, then integration by parts yields
\begin{equation}
\sum_{x\ge M} f\circ\tau_{0,x}\frac1{x^{2q}}\le c(q)f^\star \frac1{M^{2q-1}}
\end{equation}
with~$c(q)<\infty$ whenever~$2q>1$. Hence, if we assume~$q\in(1/2,1-p)$, applying this to the function $f:=\varphi\, E^0((A^\star)^{\frac q{1-p}})$ results in
\begin{equation}
\text{\rm I}_M(t)\le c(q) f^\star\, t^q \frac1{(M\sqrt t)^{2q-1}}.
\end{equation}
This shows that $\frac1{\sqrt t}\text{\rm I}_M(t)$ tends to zero as~$t\to\infty$ followed by~$M\to\infty$. The convergence occurs on~$\{f^\star<\infty\}$ which is a full-measure event because~$f\in L^1(\BbbP)$ by \eqref{E:8.7uia}.

Concerning the expression in \eqref{E:7.14}, abbreviate $t(x):=x^{2(1-p)}t^p$ and note that
\begin{equation}
P^x\bigl(A(t)\le t(x),\,Y_t<0\bigr)\le P^0\bigl(N(t(x))>2t(x)\bigr)+P^0\bigl(\max_{n\le 2t(x)}|Z_n|>x\bigr).
\end{equation}
Since~$N(t)$ is Poisson with parameter~$t$, the first probability is at most $\texte^{-c t(x)}$, for some constant~$c>0$, by a standard large-deviation estimate. The Reflection Principle in turns bounds the second probability by $2\texte^{-c x^2/t(x)}$. Bounding the sum over~$x\ge M\sqrt t$ as the sum over~$n\ge M-1$ and a sum over~$x\in[n\sqrt t,(n+1)\sqrt t)$ and invoking integration by parts shows
\begin{equation}
\begin{aligned}
\text{\rm II}_M(t)
&\le\varphi^\star
\sum_{x\ge M\sqrt t}x\bigl(\texte^{-c x^{2(1-p)}t^p}+2\texte^{-c x^{2p}/t^p}\bigr)
\\
&\le\varphi^\star\sqrt t\sum_{n\ge M-1}(n+1)\bigl(\texte^{-c n^{2(1-p)}t}+2\texte^{-c n^{2p}}\bigr)\,,
\end{aligned}
\end{equation}
where~$\varphi^\star$ is the maximal function associated with spatial shifts of~$\varphi$.
The resulting sum tends to zero as~$M\to\infty$ uniformly in~$t\ge1$.
\end{proofsect}

In order to handle the remaining part of the sums in \eqref{E:4.2e}, we will prove:

\begin{lemma}
\label{lemma-7.2}
There is~$\wh\sigma>0$ such that, for~$W:=\NN(0,\wh\sigma^2)$, $\BbbP$-a.e.\ environment and any~$M>0$,
\begin{equation}
\label{E:7.12}
\lim_{t\to\infty}\,\,\frac1{\sqrt t}\sum_{0\le x\le M\sqrt t}\varphi\circ\tau_{0,x}\,P^x(Y_t<0)
=\int_0^M\textd s\,\,P(W<-s)
\end{equation}
as well as
\begin{equation}
\label{E:7.13}
\lim_{t\to\infty}\,\,\frac1{\sqrt t}\sum_{-M\sqrt t\le x<0}\varphi\circ\tau_{0,x}\,P^x(Y_t\ge0)
=\int_{-M}^0\textd s\,\,P(W>s).
\end{equation}
\end{lemma}

Before we give the proof, note that from here we now quickly get:

\begin{proofsect}{Proof of Proposition~\ref{prop-6.2} from Lemma~\ref{lemma-7.2}}
Since the right-hand sides of \twoeqref{E:7.12}{E:7.13} coincide, Lemma~\ref{lemma-7.1} gives
\begin{equation}
\label{E:8.14o}
\lim_{t\to\infty}\frac{|\chi(-t,0)|}{\sqrt t}=0,\quad\BbbP\text{-a.s.}
\end{equation}
so we just need to turn this into a statement about the limit of times tending to positive infinity. Let~$\epsilon>0$ and set
\begin{equation}
K:=\sup_{t\ge0}\bigl(|\chi(-t,0)|-\epsilon\sqrt t\bigr).
\end{equation}
Then~$K<\infty$ $\BbbP$-a.s.\ by \eqref{E:8.14o} and so, by the Pointwise Ergodic Theorem, for $\BbbP$-a.e.\ environment there is a (random)~$R<\infty$ such that the set $\Xi_R:=\{n\in\N\colon K\circ\tau_{n,0}\le R\}$ has a positive (and well defined) density in~$\N$. This implies that there is a (random)~$n_0<\infty$ such that $\Xi_R\cap[n,2n]\ne\emptyset$ for all~$n\ge n_0$. Now assume that $t\in[n/2,n]$ for some~$n\ge n_0$ and use the above observation to find~$t_n\in\Xi_R\cap[n,2n]$. Then
\begin{equation}
\begin{aligned}
\bigl|\chi(t,0)\bigr|
&\le\bigl|\chi(t_n,0)-\chi(t,0)\bigr|+\bigl|\chi(t_n,0)\bigr|
\\
&=\bigl|\chi(t-t_n,0)\bigr|\circ\tau_{t_n,0}+\bigl|\chi(-t_n,0)\bigr|\circ\tau_{t_n,0}
\\
&\le K\circ\tau_{t_n,0}+\epsilon\sqrt{t_n-t}+K\circ\tau_{t_n,0}+\epsilon\sqrt{t_n}\,.
\end{aligned}
\end{equation}
Since $t_n\le 2n\le 4t$ and~$K\circ\tau_{t_n,0}\le R$, the right-hand side is at most $2R+2\epsilon\sqrt{4t}$. Dividing by~$\sqrt t$ and taking~$t\to\infty$ followed by~$\epsilon\downarrow0$, we  get the desired result.
\end{proofsect}

It remains to prove Lemma~\ref{lemma-7.2}. Here we will use:

\begin{lemma}
\label{lemma-7.3}
For $\BbbP$-a.e.\ realization of the environment, under $P^{0}$ we have
\begin{equation}
\frac1{\sqrt t}\,Y_{t}\,\underset{t\to\infty}{\overset{\text{\rm law}}\longrightarrow}\,\NN(0,{\wh\sigma}^2)
\end{equation}
where
\begin{equation}
\label{E:7.18}
{\wh\sigma}^2:=2\E\bigl(b_0(0)\varphi\bigr).
\end{equation}
\end{lemma}

\begin{proofsect}{Proof}
Under $P^{0}$ we have $Y_{t}= Z_{N(A(t))}$ where $t\mapsto Z_{N(t)}$ is the constant-speed conti\-nuous-time simple symmetric random walk which obeys the Functional CLT with unit limit variance. It thus suffices to show that the clock process converges to a deterministic linear function. This follows from
\begin{equation}
\frac{A(t)}{t}\,\underset{t\to\infty}{\overset{P^0}\longrightarrow}\,\wh\sigma^2,\quad \text{$\BbbP$-a.s.}
\end{equation}
which is itself proved by the Birkhoff Pointwise Ergodic Theorem applied under the stationary and ergodic law~$\Q$ and the fact that~$\Q$ is equivalent to~$\BbbP$.
\end{proofsect}

\begin{proofsect}{Proof of Lemma~\ref{lemma-7.2}}
We will again focus only on \eqref{E:7.12} as \eqref{E:7.13} is obtained analogously.
Let $\wh\sigma$ be the quantity in \eqref{E:7.18} and denote~$W:=\NN(0,\wh\sigma^2)$. Given~$\epsilon>0$, the quenched CLT for~$Y$ in Lemma~\ref{lemma-7.3} ensures there is a $\BbbP$-a.s.\ finite random variable~$T_0$ on the space of environments such that
\begin{equation}
\sup_{r\in\R}\Bigl|P^0\bigl(Y_t/\sqrt t\le r\bigr)-P(W\le r)\Bigr|<\epsilon,\quad t\ge T_0.
\end{equation}
Denote~$T_x:=T_0\circ\tau_{0,x}$ and observe that
\begin{equation}
P^x(Y_t<0)=P^0(Y_t<-x)\circ\tau_{0,x}.
\end{equation}
Decomposing the sum in \eqref{E:7.12} according to whether $\{T_x\le t\}$ occurs or not, we get
\begin{multline}
\label{E:7.22}
\qquad
\Bigl|\sum_{0\le x\le M\sqrt t}\varphi\circ\tau_{0,x}\,P^x(Y_t<0)-\sum_{0\le x\le M\sqrt t}\varphi\circ\tau_{0,x}\,P\bigl(W<-x/\sqrt t\,\bigr)\Bigr|
\\
\le\sum_{0\le x\le M\sqrt t}\varphi\circ\tau_{0,x}\bigl(\epsilon+1_{\{T_x>t\}}\bigr).
\qquad
\end{multline}
Dividing both sides by~$\sqrt t$, the Pointwise Ergodic Theorem along with the Monotone Convergence Theorem show that the right-hand side tends to zero as~$t\to\infty$ followed by~$\epsilon\downarrow0$. In light of the fact that $\chi(0,1)=\varphi-1$, Lemma~\ref{lemma-6.1} gives
\begin{equation}
\lim_{n\to\infty}\,\frac1n\sum_{x=0}^{n-1}\varphi\circ\tau_{0,x}=1,\quad\BbbP\text{-a.s.}
\end{equation}
From the monotonicity and continuity of the CDF of~$W$ we then readily get
\begin{equation}
\lim_{t\to\infty}\frac1{\sqrt t}\,\sum_{0\le x\le M\sqrt t}\varphi\circ\tau_{0,x}\,P\bigl(W<-x/\sqrt t\,\bigr)=\int_0^M \textd s\,\,P(W<-s).
\end{equation}
In combination with \eqref{E:7.22}, this now proves \eqref{E:7.12}.
\end{proofsect}

\begin{remark}
\label{rem-8.4}
Although have not quite managed to prove this, we believe that
\begin{equation}
\label{E:8.25}
\E\bigl(b_0(0)\varphi^2\bigr)=\E\bigl(b_0(0)\varphi\bigr).
\end{equation}
This is because stationarity of~$\BbbP$ under spatial shifts combined with some elementary calculus allow us to derive
\begin{equation}
\frac1t\E\bigl(\chi(t,0)^2\bigr)\,\underset{t\to\infty}\longrightarrow\,2\E\bigl(b_0(0)(\varphi-1)\varphi\bigr)
\end{equation}
and because we expect the convergence in Lemma~\ref{lemma-7.2} to hold in~$L^2(\BbbP)$-sense as well. (Alternatively, we expect $\epsilon\E(\chi_\epsilon^2)$ to vanish in the limit as~$\epsilon\downarrow0$.) If \eqref{E:8.25} indeed holds, then the limiting variance of the Brownian motion arising from the walk~$X$ is the \emph{same} as the limit variance of the Brownian motion arising from~$Y$, a fact for which we have no intuitive explanation.
\end{remark}

\section{Necessity of the moment conditions}
\label{sec8}\nopagebreak\noindent
In this final section, we will address the situations when one of the moment condition fails. We start by the lower moment condition in Theorem~\ref{lemma-lower}.
Fix~$\beta>0$ and consider the following quantity
\begin{equation}
R_\beta(t):=\frac1{t^{1/2}}\,\E\sum_{\begin{subarray}{c}
x,y\in\Z\\|x|,|y|\le\sqrt t
\end{subarray}}
\int_0^\infty\textd u\, \texte^{-\beta u}\,\,  P_a^x(X_{tu}=y)\,.
\end{equation}
The absence of the lower moment condition manifests itself as follows:

\begin{lemma}
\label{lemma-9.1}
If~$\BbbP$ is as in the statement of Theorem~\ref{lemma-lower}, then
\begin{equation}
\label{E:8.3}
\liminf_{t\to\infty} R_\beta(t)\ge\frac2\beta.
\end{equation}
\end{lemma}

\begin{proofsect}{Proof}
The proof is based on a monotonicity argument with respect to the underlying law~$\BbbP$ of static conductances.  
Let $L$ denote the generator of the random walk and let $f_t(x):=1_{[-\sqrt t,\sqrt t]}(x)$. Then
\begin{equation}
R_\beta(t)
=\frac1{t^{3/2}}\E\bigl(f_t, (\beta-L)^{-1} f_t\bigr)_{\ell^2(\Z)}.
\end{equation}
The inner product on the right-hand side is monotone decreasing with respect to the standard partial order on individual conductances and so~$R_\beta(t)$ is decreasing in~$\BbbP$. Next observe that, whenever~$\BbbP$ is such that the moment conditions in \eqref{E:1.3} hold, and~$X$ thus obeys an annealed invariance principle, we have\begin{equation}
\label{E:9.4o}
R_\beta(t)\,\underset{t\to\infty}\longrightarrow\,\int_{[-1,1]^2}\textd x\textd y\,\int_0^\infty\textd u\,\texte^{-\beta u}\,\frac1{\sqrt{2\pi\sigma^2u}}\texte^{-\frac1{2\sigma^2 u}(x-y)^2},
\end{equation}
where~$\sigma^2$ is the variance of the limiting Brownian motion. In this case~$\sigma^2$ can in fact be explicitly computed to be
\begin{equation}
\label{E:9.5o}
\sigma^2=2\bigl[\E(a(0,1)^{-1})\bigr]^{-1}
\end{equation}
thanks to the explicit representation of the corrector (see, e.g., Biskup and Prescott~\cite{BP07}).

We will now use these facts to derive the claim. Consider~$\BbbP$ as in the statement of Theorem~\ref{lemma-lower} and let~$R_\beta(t)$ be related to~$\BbbP$ as in \eqref{E:8.3}. Given~$\epsilon>0$, consider the conductance model with conductances $a^{(\epsilon)}(x,y):=a(x,y)\vee\epsilon$ and let~$R^{(\epsilon)}_\beta(t)$ be the corresponding quantity in \eqref{E:8.3}. The monotonicity in the conductance law gives
\begin{equation}
R_\beta(t)\ge R_\beta^{(\epsilon)}(t),\quad\epsilon>0.
\end{equation}
 Moreover, \twoeqref{E:9.4o}{E:9.5o} apply to~$R^{(\epsilon)}_\beta(t)$. It follows that, for any~$\epsilon>0$, the \emph{limes inferior} of~$R_\beta(t)$ is bounded from below by the right hand side of \eqref{E:9.4o} with $\sigma^2$ replaced by
\begin{equation}
\sigma_\epsilon^2:=[\E(a^{(\epsilon)}(0,1)^{-1})]^{-1}.
\end{equation}
The Monotone Convergence Theorem shows that~$\sigma_\epsilon^2\to0$ as~$\epsilon\downarrow0$ in which limit the right-hand side of \eqref{E:9.4o} tends to~$2/\beta$, as desired. 
\end{proofsect}

We are now ready to give:

\begin{proofsect}{Proof of Theorem~\ref{lemma-lower}}
In what follows we write~$\sqrt t$ instead of~$\lfloor\sqrt t\rfloor$ to ease notation. Let~$\BbbP$ be as in the statement. The key point is to prove, with the help of Lemma~\ref{lemma-9.1}, that~$R_\beta(t)$ tends to~$2/\beta$ as~$t\to\infty$. For this we first use the translation invariance of~$\BbbP$ to rewrite the desired quantity as
\begin{equation}
R_\beta(t)=\frac1{\sqrt t}\,\,\E\sum_{|x|\le 2\sqrt t+1}\bigl(2\sqrt t+1-|x|\bigr)\int_0^\infty\textd u\,\texte^{-\beta u}\, P_a^0(X_{tu}=x).
\end{equation}
If we drop the term~$|x|$, extend the sum to all~$x\in\Z$ and use that~$P_a^0(X_{tu}=\cdot)$ is a probability, we readily get
\begin{equation}
R_\beta(t)\le\frac{2\sqrt t+1}{\sqrt t}\frac1\beta.
\end{equation}
The right-hand side tends to~$2/\beta$ as~$t\to\infty$.

Lemma~\ref{lemma-9.1} now tells us that, for any~$\BbbP$ as in the statement, the approximations we used in the upper bound wash out in the limit. In particular, we must have
\begin{equation}
\frac1{\sqrt t}\,\,\E\sum_{x\in\Z}\bigl(|x|\wedge\sqrt t\,\bigr)\int_0^\infty\textd u\,\texte^{-\beta u}\,P_a^0(X_{tu}=x)
\,\underset{t\to\infty}\longrightarrow\,0.
\end{equation}
Markov's inequality readily converts this into
\begin{equation}
\label{E:8.9}
\int_0^\infty\textd u\,\,\texte^{-\beta u}\,\E P^0_a\bigl(|X_{ut}|\ge\sqrt t\bigr)\,\,\underset{t\to\infty}{\longrightarrow}\,\,0.
\end{equation}
Pick~$\delta>0$ and consider the event $\{|X_t|\ge\delta\sqrt t\}$. Let $U$ be uniform on~$[0,1]$ independent of~$a$ and~$X$ and decompose the said event according to which of the terms~$|X_{tU}|$ and~$|X_t-X_{tU}|$ is larger. A union bound combined with the Markov property for~$X$ and the invariance of~$\BbbP$ under the evolution $t\mapsto\tau_{t,X_t}(a)$ of the environment from the point of view of the particle (cf Lemma~\ref{lemma-POVP}) yield
\begin{equation}
\E P^0_a\bigl(|X_t|\ge\delta\sqrt t\bigr)
\le\int_0^1\textd u \,\,\E\Bigl[P^0\bigl(|X_{ut}|\ge\tfrac12\delta\sqrt t\bigr)+P^0\bigl(|X_{(1-u)t}|\ge\tfrac12\delta\sqrt t\bigr)\Bigr].
\end{equation}
Replacing~$t$ by~$4t/\delta^2$ now shows, via \eqref{E:8.9} and a routine change of variables, that the integral on the right-hand side tends to zero as~$t\to\infty$. The claim follows.
\end{proofsect}

Concerning the failure of the upper moment condition, we give:

\begin{proofsect}{Proof of Theorem~\ref{thm-2.2}}
Consider the spatially-homogeneous (dynamical) conductances derived from process~$\eta_t$ as in \eqref{E:2.3o}. Since the environment is homogeneous in space, the random walk~$X$ has the law of a time change of the  simple symmetric random walk. Explicitly,
\begin{equation}
X_t\,\laweq\, Z_{N(\tilde A(t))},\quad t\ge0,
\end{equation}
where $N$ is an independent rate-1 Poisson process, $Z$ is the discrete-time simple symmetric random walk on~$\Z$ and
\begin{equation}
\tilde A(t):=2\int_0^t\textd s\,\,\eta_s.
\end{equation}
The claim follows from the Central Limit Theorem for the random walk~$t\mapsto Z_{N(t)}$ and the fact that, under the assumption of ergodicity of~$t\mapsto\eta_t$ and diverging expectation of~$\eta_0$, we have $\tilde A(t)/t\to\infty$ as~$t\to\infty$ $\BbbP$-a.s.
\end{proofsect}

\section*{Acknowledgments}
\nopagebreak\nopagebreak\noindent
This project has been supported in part by the NSF award DMS-1712632 and GA\v CR project P201/16-15238S.
I am grateful to Pierre-Fran\c cois Rodriguez for valuable contributions in earlier attempts to solve this problem and, later, for keen observations on the strategy that ultimately succeeded. 
This paper is dedicated to Jean-Dominique Deuschel on the occasion of his 60th birthday. I wish to thank Jean-Dominique for his friendship and also for challenging my mathematical ability every time we meet. This note is a response to one of these challenges.

\bibliographystyle{abbrv}

\begin{thebibliography}{10}

\bibitem{A14}
S.~Andres (2014). {Invariance principle for the random conductance model with dynamic bounded conductances}. 
\textit{Ann. Inst. Henri Poincar\'e Probab. Stat.} \textbf{50}, no.~2, 352--374.

\bibitem{ACDS16}
S.~Andres, A. Chiarini, J.-D. Deuschel, M. Slowik (2018).
Quenched invariance principle for random walks with time-dependent ergodic degenerate weights.
\textit{Ann. Probab.} \textbf{46}, no.~1, 302--336


\bibitem{ASD15}
S. Andres, J.-D. Deuschel, M. Slowik (2015). 
{Invariance principle for the random conductance model in a degenerate ergodic environment}. \textit{Ann. Probab.} \textbf{43}, no.~4, 1866--1891.


\bibitem{BD12}
M.T.~Barlow and J.-D.~Deuschel (2012).
{Invariance principle for the random conductance model with unbounded conductances}. \textit{Ann. Probab.} \textbf{38} (2012), no. 1, 234--276. 

\bibitem{BB07}
N.~Berger and M.~Biskup (2007).
{Quenched invariance principle for simple random walk on percolation clusters}.
 \textit{Probab. Theory Rel. Fields} \textbf{137}, no. 1-2, 83--120.


\bibitem{B11}
M. Biskup (2011). 
{Recent progress on the Random Conductance Model}. \textit{Prob. Surveys} \textbf{8}  294--373.

\bibitem{BK-unpublished}
M. Biskup and T. Kumagai, unpublished.

\bibitem{BP07}
M.~Biskup and T.M.~Prescott  (2007).
{Functional CLT for random walk among bounded conductances}.
\textit{Electron. J. Probab.} \textbf{12}, Paper no. 49, 1323--1348.

\bibitem{BR18}
M. Biskup and P.-F. Rodriguez (2018).
Limit theory for random walks in degenerate time-dependent random environments.
\textit{J. Funct. Anal.} \textbf{274}, no. 4, 985--1046


\bibitem{DS16}
J.-D. Deuschel and M.~Slowik (2016).
{Invariance principle for the one-dimensional dynamic Random Conductance Model under moment conditions}. 
arXiv:1604.03826

\bibitem{Helland}
I. S. Helland (1982). 
Central limit theorems for martingales with discrete or continuous time. 
\textit{Scand. J. Statist.} \textbf{9}, no.~2, 79--94.

\bibitem{KV86}
C.~Kipnis and S.R.S.~Varadhan (1986). {A central limit theorem for additive functionals of reversible Markov processes and applications to simple exclusions}. \textit{Commun. Math. Phys.} \textbf{104}, no.~1, 1--19.


\bibitem{Kumagai}
T.~Kumagai (2014). \textit{Random walks on disordered media and their scaling limits}, Lecture notes from the 40th Probability Summer School held in Saint-Flour, 2010. Lecture Notes in Mathematics, 2101. \'Ecole d'\'Et\'e de Probabilit\'es de Saint-Flour. Springer, Cham, x+147 pp.

\bibitem{Liggett-MC}
T.M. Liggett (2010). \textit{Continuous time Markov processes: An introduction}, Graduate Studies in Mathematics, vol. 113, Amer. Math. Soc., 271 pp. 

\end{thebibliography}

\end{document}